\nonstopmode
\documentclass[12pt]{amsart}
\usepackage{amscd}
\usepackage{amssymb}

\numberwithin{equation}{section}

\swapnumbers

\theoremstyle{plain} 
\newtheorem{thm}[equation]{Theorem}

\newtheorem{lem}[equation]{Lemma}
\newtheorem{prop}[equation]{Proposition}

\theoremstyle{definition}
\newtheorem{defn}[equation]{Definition}

\theoremstyle{remark}
\newtheorem{rem}[equation]{Remark}

\newcommand{\Observation}{General principle}

\newtheorem*{VariableNoNum}{{\VariableText}}
\newtheorem{Variable}[equation]{{\VariableText}}

\theoremstyle{definition}
\newtheorem*{VariableNoNumBold}{{\VariableText}}
\newtheorem{VariableBold}[equation]{{\VariableText}}

\newenvironment{titled}[1]
     {\def\VariableText{{#1}}\begin{VariableNoNum}}
     {\end{VariableNoNum}}

\newenvironment{numbered}[1]
     {\def\VariableText{{#1}}\begin{Variable}}
     {\end{Variable}}

\newenvironment{NumberedSubSection}[1]
     {\def\VariableText{{\textrm{\textbf{#1}}}}\begin{VariableBold}}
     {\end{VariableBold}}

\newlength{\asidelength}

\def\Changed/{\ifvmode\else\vadjust{%
\vbox to 0pt{\vskip -\baselineskip%
\hbox to 0pt{\hss\vrule height 0pt depth 1.2\baselineskip\hskip 1em}\vss}}\fi}
\def\CHanged{\ifvmode\else\vadjust{%
\vbox to 0pt{\vskip -\baselineskip%
\hbox to 0pt{\hss\vrule height 0pt depth 1.2\baselineskip\hskip 1em}\vss}}\fi}

\def\Math#1{\def\MathString{#1}\futurelet\MathDelim\MathChoose}
\def\MathChoose{\ifmmode\let\MathDo\MathString%
              \else\let\MathDo\MathSkip\fi%
              \MathDo}
\def\MathSkip{\ifx\MathDelim/\def\MathDo{$\MathString$\EatOne}%
              \else\def\MathDo{$\MathString$}\fi%
              \MathDo}

\def\Text#1{\def\TextString{#1}\futurelet\TextDelim\TextSkip}
\def\TextSkip{\ifx\TextDelim/\def\TextDo{\TextString\EatOne}%
              \else\let\TextDo\TextString\fi%
              \TextDo}
\def\EatOne#1{}

\def\SkipToEndScan#1\EndScan{}
\def\Scan#1#2#3{\ifx#1#2#3\expandafter\SkipToEndScan\fi\Scan#1}
\def\Upper#1{%
\Scan#1aAbBcCdDeEfFgGhHiIjJkKlLmMnNoOpPqQrRsStTuUvVwWxXyYzZ#1#1\EndScan}
\def\Phrase#1 #2/#3/#4=#5 #6/#7/#8.{%
\expandafter\edef\csname#2#3\endcsname{\noexpand\Text{#6#7}}
\expandafter\edef\csname\Upper#2#3\endcsname{\noexpand\Text{\Upper#6#7}}
\expandafter\edef\csname#1#2#3\endcsname{\noexpand\Text{#5 #6#7}}
\expandafter\edef\csname\Upper#1#2#3\endcsname{\noexpand\Text{\Upper#5 #6#7}}
\expandafter\edef\csname#2#4\endcsname{\noexpand\Text{#6#8}}
\expandafter\edef\csname\Upper#2#4\endcsname{\noexpand\Text{\Upper#6#8}}
}

\newcommand{\colim}{\operatorname{colim}}
\newcommand{\hocolim}{\operatorname{hocolim}}
\newcommand{\holim}{\operatorname{holim}}

\newcommand{\whatever}{\text{--}}

\newcommand{\op}{^{\text{op}}}

\newcommand{\Z}{\Math{\mathbb{Z}}}
\newcommand{\Q}{\Math{\mathbb{Q}}}

\newcommand{\Fp}{\Math{\mathbb{F}_p}}
\newcommand{\Zpinfty}{\Math{\Z/p^\infty}}
\newcommand{\Sphere}{\Math{\mathbb{S}}}
\newcommand{\Zpadic}{\Math{\mathbb{Z}_p}}

\newcommand{\Smash}{\wedge}
\newcommand{\Tensor}{\otimes}

\newcommand{\RightArrow}[1]{\xrightarrow{#1}} 

\newcommand{\Hom}{\operatorname{Hom}}

\newcommand{\Map}{\operatorname{Map}}
\newcommand{\End}{\operatorname{End}}
\newcommand{\Ext}{\operatorname{Ext}}
\newcommand{\Tor}{\operatorname{Tor}}

\newcommand{\hhh}{\operatorname{h}\!}
\newcommand{\thhh}{\tilde{\operatorname{h}}}
\def\HomotopyOrbit#1on#2/{\ensuremath{#2_{\hhh#1}}}
\def\RedHomotopyOrbit#1on#2/{\ensuremath{#2_{\thhh#1}}}

\newcommand{\iso}{\cong}
\newcommand{\weq}{\sim}

\def\NOTE/{}
\def\Phrase#1 #2/#3/#4=#5 #6/#7/#8.{%
\expandafter\edef\csname#2#3\endcsname{\noexpand\Text{#6#7\noexpand\NOTE/}}
\expandafter\edef\csname\Upper#2#3\endcsname{\noexpand\Text{\Upper#6#7\noexpand\NOTE/}}
\expandafter\edef\csname#1#2#3\endcsname{\noexpand\Text{#5 #6#7\noexpand\NOTE/}}
\expandafter\edef\csname\Upper#1#2#3\endcsname{\noexpand\Text{\Upper#5 #6#7\noexpand\NOTE/}}
\expandafter\edef\csname#2#4\endcsname{\noexpand\Text{#6#8\noexpand\NOTE/}}
\expandafter\edef\csname\Upper#2#4\endcsname{\noexpand\Text{\Upper#6#8\noexpand\NOTE/}}
}

\Phrase an ehregular//s=a proxy-small//s.
\newcommand{\ehregularity}{\Text{proxy-smallness\NOTE/}}
\Phrase a regular//s=a small//s.
\Phrase a cohregular//s=a cosmall//s.
\Phrase an hregular//s=a small//s.
\Phrase a mixedregularit/y/ies=a proxy-smallness//es.

\def\Phrase#1 #2/#3/#4=#5 #6/#7/#8.{%
\expandafter\edef\csname#2#3\endcsname{\noexpand\Text{#6#7}}
\expandafter\edef\csname\Upper#2#3\endcsname{\noexpand\Text{\Upper#6#7}}
\expandafter\edef\csname#1#2#3\endcsname{\noexpand\Text{#5 #6#7}}
\expandafter\edef\csname\Upper#1#2#3\endcsname{\noexpand\Text{\Upper#5 #6#7}}
\expandafter\edef\csname#2#4\endcsname{\noexpand\Text{#6#8}}
\expandafter\edef\csname\Upper#2#4\endcsname{\noexpand\Text{\Upper#6#8}}
}

\newcommand{\Salgebra}{$\mathbb{S}$-algebra}
\Phrase a ring//s=an \relax\relax\noexpand\Salgebra//s.
\newcommand{\dccomplete}{\Text{dc--complete}}
\newcommand{\dccompleteness}{\Text{dc-completeness}}
\Phrase an hGorenstein//s=a Gorenstein//s.
\Phrase an augmented//s=an augmented//s.
\Phrase a smalll//s=a small//s.
\Phrase a cosmall//s=a cosmall//s.
\Phrase a mixedsmall//s=a proxy-small//s.
\Phrase a lift//s=a lift//s.
\newcommand{\Elift}{\Text{$\EE$-lift}}
\newcommand{\Elifts}{\Text{$\EE$-lifts}}

\newcommand{\DGAs}{\Text{DGAs}}

\newcommand{\Derived}{\mathbf{D}}

\newcommand{\HomOver}[1]{\Hom_{#1}}
\newcommand{\EndOver}[1]{\End_{#1}}
\newcommand{\HomR}{\HomOver R}
\newcommand{\HomE}{\HomOver{\EE}}
\newcommand{\Homk}{\HomOver k}
\newcommand{\HomS}{\HomOver\Sphere}

\newcommand{\EndR}{\EndOver R}

\newcommand{\TensorOver}[1]{\Tensor_{#1}}
\newcommand{\TensorR}{\TensorOver R}
\newcommand{\TensorE}{\TensorOver {\EE}}
\newcommand{\TensorS}{\TensorOver{\Sphere}}

\newcommand{\Tensork}{\TensorOver{k}}
\newcommand{\EndE}{\EndOver{\EE}}

\newcommand{\RIhat}{R\hat{{}_I}}

\newcommand{\Suspend}{\Sigma}
\newcommand{\Stable}{\Suspend^\infty}
\newcommand{\Shift}{\Sigma}

\newcommand{\SH}{a}                      

\newcommand{\E}{E}

\newcommand{\EE}{\Math{\mathcal{E}}}
\newcommand{\EEK}{\Math{\EE_K}}

\newcommand{\T}{T}

\newcommand{\Pair}[2]{(#1,#2)}
\newcommand{\DCellOfPair}[2]{\Derived\!\Cell\Pair{#1}{#2}}

\newcommand{\Cell}{\operatorname{Cell}}
\newcommand{\CellOf}[1]{\Cell_{#1}}
\newcommand{\CellOverOf}[2]{\Cell^{#1}_{#2}}
\newcommand{\Cellk}{\CellOf{k}}

\newcommand{\II}{\Math{\mathcal{I}}}

\newcommand{\Hof}[1]{\pi_{#1}}

\newcommand{\Modules}{\operatorname{Mod}}
\newcommand{\OfMod}[1]{{}_{#1}\Modules}
\newcommand{\ModOf}[1]{\Modules_{#1}}
\newcommand{\RMod}{\Math{\OfMod{R}}}
\newcommand{\ModR}{\Math{\ModOf{R}}}

\newcommand{\Ho}{\operatorname{Ho}}

\newcommand{\Boxtensor}{\boxtimes}
\newcommand{\aBoxtensor}{\boxtimes^a}
\begin{document}
\title[Duality]{Duality in algebra and topology}

\author{W. G. Dwyer}
\author{J. P. C. Greenlees}
\author{S. Iyengar}

\address{Department of Mathematics, University of Notre Dame, Notre
  Dame, IN 46556. USA}
\address{Department of Pure Mathematics, Hicks Building, Sheffield S3
  7RH. UK}
\address{Department of Mathematics, University of Nebraska,
  Lincoln, NE 68588. USA}

\email{dwyer.1@nd.edu}
\email{j.greenlees@sheffield.ac.uk}
\email{iyengar@math.unl.edu}

\date{\today}

\dedicatory{Dedicated to Clarence W. Wilkerson, on the occasion of his
  sixtieth birthday}

\maketitle

\tableofcontents

\section{Introduction}\label{SIntroduction}

In this paper we take some classical ideas from commutative algebra,
mostly ideas involving duality, and apply them in algebraic topology.
To accomplish this we interpret properties of ordinary commutative
rings in such a way that they can be extended to the more general
rings that come up in homotopy theory. Amongst the rings we work with
are the differential graded ring of cochains on a space $X$, the
differential graded ring of chains on the loop space $\Omega X$, and
various ring spectra, e.g., the Spanier-Whitehead duals of finite
spectra or chromatic localizations of the sphere
spectrum.

Maybe the 
most important contribution of this paper is the conceptual framework,
which allows us to view all of the following dualities
\begin{itemize}
\item Poincar\'e duality for manifolds
\item Gorenstein duality for commutative rings
\item Benson-Carlson duality for cohomology rings of finite groups
\item Poincar\'e duality for groups
\item Gross-Hopkins duality in chromatic stable homotopy theory
\end{itemize}
as examples of a \emph{single} phenomenon. Beyond setting up this
framework, though, we prove some new results, both in algebra and
topology, and give new proofs of a number of old results. Some of the
rings we look at, such as $C_*\Omega X$, are not commutative in any
sense, and so implicitly we extend the methods of commutative algebra
to certain noncommutative settings.  We give a new formula for the
dualizing module of a Gorenstein ring (\ref{LocalRingMatlisLift}); this formula involves
differential graded algebras (or ring spectra) in an essential way
and is one instance of a general construction that in another setting
gives the Brown-Comenetz dual of the sphere spectrum (\ref{BrownComenetz}).  We also prove
the local cohomology theorem for $p$-compact groups
and reprove it for compact Lie groups with orientable adjoint representation (\ref{GoodLoopSpace}).
The previous proof for compact Lie groups
\cite{benson-greenlees;commutative} uses equivariant topology, but ours
does not.

\begin{NumberedSubSection}{Description of results}
  The objects we work with are fairly general; briefly, we allow
  rings, differential graded algebras (DGAs), or ring spectra; these
  are all covered under the general designation \emph{\ring} (see
  \ref{OurNotation}). We usually work in a derived category or
  in a homotopy category of module spectra, to the extent that even if
  $R$~is a ring, by a \emph{module} over~$R$ we mean a chain complex of ordinary $R$-modules. Most of the time we start
  with a homomorphism $R\to k$ of \rings/ and let $\EE$ denote the
  endomorphism \ring/ $\EndR(k)$. There are three main parts to the
  paper, which deal with three different but related types of
  structures: \emph{smallness}, \emph{duality}, and the \emph{Gorenstein condition}.

  \begin{titled}{Smallness} There are several different kinds of
    smallness which the homomorphism $R\to k$ might enjoy
    (\ref{DefineRegular});  the
    weakest and most flexible one is called \emph{\mixedregularity/}. Any
    surjection from a commutative Noetherian ring to a regular ring is \ehregular/
    (\ref{OrdinaryRingERegular}).  One
    property of \anehregular/ homomorphism is particularly interesting
    to us.  Given an $R$-module $M$, there is an associated module
    $\Cellk(M)$, which is the closest $R$-module approximation to $M$
    which can be cobbled together from shifted copies of $k$ by using
    sums and exact triangles.     The notation $\Cellk(M)$ comes from topology \cite{Farjoun}, but
    if $R$ is a commutative ring and $k=R/I$ for a finitely generated
    ideal $I$, then $\pi_*\Cellk(M)$ is just the local
    cohomology $H^{-*}_I(M)$ of $M$ at $I$ \cite[\S6]{dwyer-greenlees}.
It turns out that if $R\to k$ is \ehregular/, there is a
    canonical equivalence (\ref{TLatentMoritaEquiv})
    \begin{equation}\label{BasicCellConstruction}
        \CellOf k M\weq \HomR(k,M)\TensorE k\,.
    \end{equation}
  \end{titled}

  \begin{titled}{Duality}
     Given $R\to k$, we look for a notion of ``Pontriagin duality'' over
     $R$ which extends the notion of ordinary duality over $k$; more
     specifically, we look for an $R$-module $\II$ such that
     $\Cellk(\II)\weq\II$ and such that for any
     $k$-module $X$, there is a natural weak equivalence
     \begin{equation}\label{BasicPontriaginEquation}
         \HomR(X,\II) \weq \Homk(X,k)\,.
     \end{equation}
     The associated Pontriagin duality (or \emph{Matlis duality}) for
     $R$-modules sends $M$ to $\HomR(M,\II)$. If $R\to k$ is
     $\Z\to\Fp$, there is only one such $\II$, namely $\Zpinfty(=\Z[1/p]/\Z)$, and
     $\Hom_{\Z}(\whatever,\Zpinfty)$ is ordinary $p$-local Pontriagin
     duality for abelian groups. Guided by a combination of
     \ref{BasicCellConstruction} and \ref{BasicPontriaginEquation}, we find that in many circumstances,
     and in particular if $R\to k$ is \ehregular, such dualizing
     modules $\II$ are determined by \emph{right} $\EE$-module
     structures on $k$; this structure is a new bit of information,
     since in its state of nature $\EE$ acts on $k$ from the left.
     Given a suitable right action, the dualizing module $\II$ is
     given by the formula
    \begin{equation}\label{BasicDualizingModuleFormula}
        \II \weq k\TensorE k\,\,,
    \end{equation}
    which mixes the exceptional right action of $\EE$ on $k$ with the
    canonical left action. This is a formula which in one setting
    constructs the injective hull of the residue class field
    of a local ring (\ref{LocalRingMatlisLift}), and in another gives the $p$-primary component
    of the Brown-Comenetz dual of the sphere spectrum (\ref{BrownComenetz}). There are also
    other examples (\S\ref{MatlisLiftingExamples}). We call an $\II$
    which is of the form described in \ref{BasicDualizingModuleFormula} a \emph{Matlis lift} of~$k$.
  \end{titled}

  \begin{titled}{The Gorenstein condition}
    The homomorphism $R\to k$ is said to be \emph{\hGorenstein} if, up
    to a shift, $\Cellk(R)$ is a Matlis lift of~$k$. This amounts to
    requiring that
    $\HomR(k,R)$ be equivalent to a shifted module $\Shift^\SH k$, and that the right action of $\EE$ on $k$ provided by this
    equivalence act as in  \ref{BasicDualizingModuleFormula} to give  a
    dualizing module~$\II$.
    There are several consequences of the \hGorenstein/ condition. 
    In the commutative ring case with $k=R/I$, 
    the equivalences
    \[
      \II= k\TensorE k =\Shift^{-\SH}\HomR(k,R)\TensorE k\weq\Shift^{-\SH}\CellOf
      kR\,.
    \]
    give a connection between the
    dualizing module $\II$ and the local cohomology object
    $\mathbf{R}\Gamma_I(R)=\CellOf    kR$.  (Another notation for
    $\mathbf{R}\Gamma_I(R)$ might be $H^\bullet_I(R)$, since $\pi_i\mathbf{R}\Gamma_I(R)$
    is the local cohomology group
    $H^{-i}_I(R)$.) In this paper we head in a slightly different direction.
    Suppose that $R$ is \anaugmented/ $k$-algebra and $R\to k$ is the
    augmentation; in this case it is possible to compare the two right
    $\EE$-modules $\HomR(k,R)$ and $\HomR(k,\Homk(R,k))$. Given that
    $R\to k$ is \hGorenstein, the first is abstractly equivalent to
    $\Shift^ak$; the second, by an adjointness argument, is always
    equivalent to $k$. If these two objects are the same as
    $\EE$-modules after the appropriate shift, we obtain a formula
    \[
           \Shift^\SH\CellOf k\Homk(R,k) \weq \CellOf k R\,,
    \]
    relating duality on the left to local cohomology on the right. 
    In many circumstances $\CellOf k\Homk(R,k)$ is equivalent to
    $\Homk(R,k)$ itself, and in these cases the above formula becomes
    \[
          \Shift^\SH\Homk(R,k)\weq \CellOf kR\,.
    \]
    This leads to spectral sequences relating the local cohomology of
    a ring to some kind of $k$-dual of the ring, for instance, if $X$
    is a suitable space, relating the local cohomology of $H^*(X;k)$
    to $H_*(X;k)$. We use this approach to reprove the local
    cohomology theorem for compact Lie groups and prove it for
    $p$-compact groups.
  \end{titled}
    
   \nocite{dgi:algebra}\nocite{dgi:chromatic}
   We intend to treat the special case of
   chromatic stable homotopy theory in
   \cite{dgi:chromatic}; it turns out that Gross-Hopkins duality is a
   consequence of the fact that the \ring/ map $L_{K(n)}\Sphere\to
   K(n)$ is Gorenstein. In
   \cite{dgi:algebra} we use our techniques to study derived
   categories of local rings.
\end{NumberedSubSection}


  \begin{NumberedSubSection}{Notation and terminology}\label{OurNotation}
    In this paper we use the term \emph{\ring/} to mean ring spectrum in
    the sense of \cite{EKMM2} or \cite{HSS}; the symbol \Sphere/
    stands for the sphere spectrum. If $k$ is a commutative \ring/, we
    refer to algebra spectra over $k$ as \emph{$k$-algebras}. The
    sphere \Sphere/ is itself a commutative ring spectrum, and, as the
    terminology ``\ring/'' suggests, any ring spectrum is 
    an algebra spectrum over \Sphere. There is a brief introduction to
    the machinery of \rings/ in \S\ref{CSpectra}; this follows the
    approach of \cite{HSS}.

    Any ring $R$ gives rise to \aring/ (whose homotopy is 
    $R$, concentrated in degree~$0$), and we do not make a distinction
    in notation between $R$ and this associated spectrum.  If $R$ is
    commutative in the usual sense it is also commutative as \aring;
    the category of $R$-algebras (in the way in which we use the term)
    is then equivalent to the more familiar category of differential
    graded algebras (DGAs) over $R$. For instance, \Z-algebras are
    essentially DGAs; \Q-algebras are DGAs over the rationals.  
See \cite{ShipleyDGA} for a
    detailed treatment of the relationship between \Z-algebras and
    DGAs.

A
    \emph{module} $M$ over \aring/ $R$ is for us a module spectrum
    over $R$; the category of these is denoted \RMod.
   Note that
      unspecified modules are left modules. 
If $R$ is a ring, then an $R$-module in our sense is essentially
    an unbounded chain complex over $R$. More generally, if $R$ is a
    \Z-algebra, an $R$-module is essentially the same as a
    differential graded module over the corresponding DGA \cite{ShipleyDGA}.
     (Unbounded chain
    complexes over a ring should be treated homologically as in
    \cite{Spaltenstein}. Differential graded modules over
    a DGA are treated very similarly; there are implicit discussions
    of this in
    \cite[\S3]{ShipleySchwedeEquiv} and \cite[\S2]{ShipleyDGA}.)  If
    $R$ is a ring, 
    any ordinary module $M$ over $R$
    gives rise to an $R$-module in our sense by the analog of the usual device
    of treating $M$ as a chain complex concentrated in
    degree~$0$. We will refer to such an $M$ as a \emph{discrete}
    module over $R$, and we will not distinguish in notation between
    $M$ and its associated spectrum. 

    \begin{titled}{Homotopy/homology}
      The homotopy groups of \aring/ $R$ and an $R$-module $M$ are
      denoted respectively $\Hof*R$ and $\Hof*M$. The group $\Hof 0R$
      is always a ring, and a ring is distinguished among
      \rings/ by the fact that $\Hof iR\iso 0$ for $i\ne0$. If $R$ is
      a \Z-algebra and $M$ is an $R$-module, the homotopy groups
      $\Hof*R$ and $\Hof*M$ amount to the homology groups of the
      corresponding differential graded objects.  A homomorphism $R\to
      S$ of \rings/ or $M\to N$ of modules is an \emph{equivalence}
      (weak equivalence, quasi-isomorphism) if it induces an
      isomorphism on $\Hof*$. In this case we write $R\weq S$ or
      $M\weq N$.  
    \end{titled}

    \begin{titled}{Hom and tensor}
      Associated to two $R$-modules $M$ and $N$ is a spectrum
      $\HomR(M,N)$ of homomorphisms; each $R$-module $M$ also has an
      endomorphism ring $\EndR(M)$. These are derived objects; for
      instance, in forming $\EndR(M)$ we always tacitly assume that
      $M$ has been replaced by an equivalent $R$-module which is
      cofibrant (projective) in the appropriate sense.   If $M$ and $N$ are
      respectively right and left modules over $R$, there is a derived
      smash product, which corresponds to tensor product of
      differential graded modules, and which we write $M\TensorR N$.

      To fix ideas, suppose that $R$ is a ring, $M$ is a discrete
      right module over $R$, and $N, K$ are discrete left modules. Then
      $\Hof i(M\TensorR N)\iso\Tor_i^R(M,N)$, while $\Hof i\HomR(K,N)\iso\Ext^{-i}_R(K,N)$. In this situation we sometimes
      write $\hom_R(M,N)$ (with a lower-case ``h'') for the group
      $\Ext^0_R(M,N)$ of ordinary $R$-maps $M\to N$.

      There are other contexts as above in which we follow the practice of
      tacitly replacing one object by an equivalent one without
      changing the notation. For instance, suppose that $R\to k$ is a
      map of \rings, and let $\EE=\EndR(k)$. The right action of $k$
      on itself commutes with the left action of $R$, and so produces
      what we refer to as a ``homomorphism $k\op\to \EE$'', although
      in general this homomorphism can be realized as a map of \rings/
      only after adjusting $k$ up to weak equivalence. The issue is
      that in order to form $\EndR(k)$, it is necessary to work with a
      cofibrant (projective) surrogate for $k$ as a left $R$-module,
      and the right action of $k$ on itself cannot in general be
      extended to an action of $k$ on such a surrogate without tweaking
      $k$ to some extent. The reader might want to consider
      the example $R=\Z$, $k=\Fp$ from \cite[\S3]{dwyer-greenlees},
      where it is clear that the ring \Fp/ cannot map to the DGA
      representing $\EE$, although a DGA weakly equivalent to $\Fp$
      does map to $\EE$. In general we silently pass over these
      adjustments and replacements in order to keep the exposition
      within understandable bounds.
    \end{titled}

    \begin{titled}{Derived category}
      The \emph{derived category} $\Derived(R)=\Ho(\RMod)$ of \aring/ $R$ is
      obtained from \RMod/ by formally inverting the weak
      equivalences.  A map between $R$-modules passes to an
      isomorphism in $\Derived(R)$ if and only if it is a weak
      equivalence. Sometimes we have to consider a homotopy category
      $\Ho(\ModR)$ involving right $R$-modules; since a right
      $R$-module is the same as a left module over the opposite ring
      $R\op$, we write $\Ho(\ModR)$ as $\Derived(R\op)$. If $R$ is a
      ring, $D(R)$ is categorically equivalent to the usual derived
      category of $R$.
    \end{titled}

    \begin{titled}{Augmentations} Many of the objects we work with are
      augmented. An \emph{augmented} $k$-algebra $R$ is a $k$-algebra
      together with an augmentation homomorphism $R\to k$ which splits
      the $k$-algebra structure map $k\to R$. A map of augmented
      $k$-algebras is a map of $k$-algebras which respects the
      augmentations. If $R$ is an augmented $k$-algebra, we will
      by default treat $k$ as an $R$-module via the homomorphism $R\to
      k$.
    \end{titled}
  
    \begin{titled}{Another path}
      The advantage of using the term \ring/ is that we can refer to
      rings, DGAs, and ring spectra in one breath. The reader
      can confidently take $\Sphere=\Z$, read DGA for \ring, $H_*$ for
      $\Hof*$, and work as in \cite{dwyer-greenlees} in the algebraic
      context of \cite{Spaltenstein}; only some examples will be lost.
      Note however that the loss will include all examples involving
      commutativity in any essential way, unless the commutative \rings/ in question are
      \Q-algebras or ordinary commutative rings. This is a consequence of the fact that
      under the correspondence between \Z-algebras and DGAs, the notion
      of commutativity for \Z-algebras does not carry over to the
      usual notion of commutativity for DGAs,
      except in characteristic~$0$ \cite[App.~C]{MandellPadic}
      \cite{kriz-may}, or if the homotopy of the \Z-algebra is
      concentrated in degree~$0$. 
    \end{titled}

  \end{NumberedSubSection}


  \begin{NumberedSubSection}{Organization of the paper}
    Section \ref{CSpectra} has a brief expository introduction to
    spectra and \rings, and section \ref{CModuleProperties} describes
    some elementary properties of \rings/ which we use later on.
    Some readers may wish to skip these sections the first time
    through. 
    The three main themes, smallness, duality, and the Gorenstein
    condition, are treated respectively in Sections \ref{CContext},
    \ref{CMatlisDuality}, and \ref{CGorenstein}.  Section
    \ref{LocalCohomologyTheorem} explains how to set up a local
    cohomology spectral sequence for a suitable \hGorenstein/
    \ring. We spend a lot of time dealing with examples;
    \S\ref{CRegularityExamples} has examples relating to smallness,
    \S\ref{MatlisLiftingExamples} examples related to duality, and
    \S\ref{GorensteinExamples} examples related to the Gorenstein
    condition. In particular, Section \ref{GorensteinExamples}
    contains a proof of 
    the local cohomology theorem for $p$-compact groups
    (\ref{GoodLoopSpace}) and for compact Lie groups with orientable
    adjoint representation
    (\ref{GorensteinCompactLie}); following
    \cite{rGreenleesCommutative} and
    \cite{benson-greenlees;commutative}, for finite groups this is one version of
    Benson-Carlson duality \cite{benson-carlson}. 
  \end{NumberedSubSection}

\begin{NumberedSubSection}{Relationship to previous work}
  There is a substantial literature on Gorenstein rings. Our
  definition of \anhGorenstein/ map $R\to k$ of \rings/ extends the
  definition of Avramov-Foxby \cite{avramov-foxby} (see
  \ref{SpecialGorensteinConditions}). F\'elix, Halperin, and Thomas
  have considered pretty much this same extension in the topological
  context of rational homotopy theory and DGAs \cite{halperin-thomas};
  we generalize their work and have benefitted from it. Frankild and
  Jorgensen \cite{frankild} have also studied an extension of the
  Gorenstein condition to DGAs, but their intentions are quite
  different from ours.
\end{NumberedSubSection}

\section{Spectra, $\mathbb S$-algebras, and commutative $\mathbb S$-algebras\relax}\label{CSpectra}

In this section we give a brief introduction to the nuts and bolts of
spectra, \rings, and commutative \rings; this is purely expository,
and so the experts can safely move on to \S\ref{CModuleProperties}.  In very rough terms, a
spectrum is something like an unbounded chain complex in which the
commutative and associative laws for the addition of elements hold
only up to coherently specified homotopies. As will become clear
later, a spectrum could also reasonably be styled an \Sphere-module.
There is a well--behaved tensor product (AKA smash product) for
spectra, and with the help of this it is possible to give simple
definitions of \rings/ and commutative \rings.

Spectra are based in one way or another on homotopy--theoretic
objects; the specific homotopy--theoretic objects we pick are
\emph{pointed simplicial sets} \cite{rGoerssJardine}. For the rest of this section,
the word \emph{space} taken by itself will mean pointed simplicial
set.  (In the course of the paper we often refer to topological
spaces, but there are standard constructions
which make it possible to pass
back and forth between topological spaces and simplicial sets without
losing homotopical information.) To
set up the category of spectra we will rely on the \emph{symmetric
  spectrum} machinery of Hovey, Shipley and Smith \cite{HSS}. Both
simplicial sets and symmetric spectra are inherently combinatorial
objects, and so from the point of view we are taking a spectrum is
combinatorial, or even algebraic, in nature.

\begin{NumberedSubSection}{Asymmetric spectra}

  We will start with a simple construction, which long predates the
  notion of symmetric spectrum; for want of a better term we'll call
  the objects that come up ``asymmetric spectra''. The category of
  asymmetric spectra is a good homotopy theoretic model for the
  category of spectra, but it lacks a decent tensor product; passing
  to the more complicated category of symmetric spectra will
  solve this problem.

  \begin{defn}\label{AsymmetricBoxtensor}
    A \emph{sequence} $X$ of spaces is a collection
    $\{X_n\}_{n\ge0}$ of  spaces. The (graded) product
    $X\aBoxtensor Y$ of two such sequences is the sequence $Z$ given by
    \[           Z_n= \coprod_{i+j=n} X_i\Smash Y_j\,.    \]
  \end{defn}

  Here $\Smash$ denotes the smash product of spaces. The superscript
  $a$ in $\aBoxtensor$ signifies ``asymmetric''; later on we will
  define another kind of $\Boxtensor$. In spite of its decoration, the
  operation $\aBoxtensor$ gives a symmetric monoidal structure on the
  category of sequences of spaces; the unit is the sequence
  $\epsilon$ with $\epsilon_0=S^0$ and $\epsilon_i=*$ for $i>0$ (here
  $S^0$ is the zero-sphere, i.e., the unit for the smash product
  operation on the category of spaces).  The twist isomorphism
  $X\aBoxtensor Y\iso Y\aBoxtensor X$ acts at level $n$ by using the
  usual isomorphisms
  $X_i\Smash Y_j\to Y_j\Smash X_i$ ($i+j=n$). Let $S^1$ denote the
  simplicial circle, and $S^n$ $(n\ge1)$ the smash power
  $S^1\Smash\cdots \Smash S^1$ ($n$~times).  There is a sequence $S$
  with $S_i=S^i$, and it is easy to produce a pairing map
  \[
              S\aBoxtensor S\to S
  \]
  which makes $S$ into a monoid with respect to $\aBoxtensor$
  \cite[2.3.4]{HSS}. This pairing map is constituted from the obvious
  isomorphisms $S^i\Smash S^j\to S^{i+j}$.

  \begin{defn}\label{DefineAsymmetricSpectrum}
    An \emph{asymmetric spectrum} is a sequence $X$ of spaces
    together with a left action of $S$ on $X$, i.e, a map $S\aBoxtensor
    X\to X$ which satisfies appropriate associativity and unital
    identities.
  \end{defn}

  \begin{rem}
    Since $S$ is the free monoid with respect to $\aBoxtensor$ on a copy
    of $S^1$ at level~$1$, an asymmetric spectrum amounts to a
    sequence $X$ of  spaces together with structure maps $S^1\Smash
    X_i\to X_{i+1}$. This is exactly a spectrum in the sense of
    Whitehead \cite{rGWWhitehead} (although he worked with topological
    spaces instead of with simplicial sets).
  \end{rem}

  If $X$ is an asymmetric spectrum, the structure maps $S^1\Smash X_i\to X_{i+1}$ induce
  homotopy group maps $\pi_nX_i\to\pi_{n+1}X_{i+1}$.

  \begin{defn}\label{ASweakequivalence}
    The \emph{homotopy groups} of an asymmetric spectrum $X$ are the
    groups $\pi_nX=\colim_i\pi_{n+i}X_i$. A map $X\to Y$ of asymmetric
    spectra is a \emph{weak equivalence} (quasi--isomorphism) if it
    induces isomorphisms $\pi_iX\iso\pi_iY$, $i\in\Z$.
  \end{defn}

  \begin{titled}{Additivity, associativity, commutativity}
    There is an elaborate homotopy theory of asymmetric spectra based
    upon this definition of weak equivalence. On its own, the
    definition of weak equivalence subtly imposes the additivity,
    associativity, and commutativity structures referred to in the
    introduction to the section.  For instance, it follows from
    \ref{ASweakequivalence} that any asymmetric spectrum $X$ is weakly
    equivalent to $\Omega^n ( S^n\Smash X)$ (where the loop functor
    $\Omega^n$ and the smash functor $S^n\Smash\text{--}$ are applied
    levelwise).  However, for any space $A$, $\Omega^n(S^n\Smash A)$ ($n\ge1$)
    has up to homotopy an associative multiplication, and these
    multiplications  enjoy ever richer
    commutativity properties as $n$~increases.
  \end{titled}

  \begin{numbered}{Relationship to chain complexes}\label{AsymmetricChains}
  As defined in \ref{ASweakequivalence}, asymmetric spectra have both positive
  and negative dimensional homotopy groups, just as unbounded chain
  complexes have both positive and negative dimensional homology
  groups. An unbounded chain complex $C$ can be converted
  into an asymmetric spectrum $X$ by setting
  $X_i=N^{-1}(\operatorname{t}\Sigma^iC)$, where $\Sigma^i$ shifts the
  complex upward $i$~times, ``$\operatorname{t}$'' cuts off the
  negative dimensional components, and $N^{-1}$ is the Dold-Kan
  denormalization functor, which converts a nonnegative chain complex
  into a simplicial abelian group.  Then $\pi_iX\iso H_iC$.
  \end{numbered}

  \begin{numbered}{Homotopy category, shifting, triangulated
      structure, suspension spectra}\label{AsymmetricShift}
    The \emph{homotopy category} of asymmetric spectra is constructed
    by formally inverting weak equivalences, just as the derived
    category of a ring is constructed by formally inverting
    quasi--isomorphisms between chain complexes.
    For any $i$ there is a shift map $\Sigma^i$ defined on the
    category of asymmetric spectra, given by $(\Sigma^iX)_n = X_{n-i}$;
    the formula is to be interpreted to mean that if $n-i<0$, then
    $(\Sigma^iX)_n=*$. This is parallel to the usual shift operation
    on chain complexes: in particular, $\pi_k\Sigma^iX=\pi_{k-i}X$,
    and up to weak equivalence $\Sigma^i\Sigma^j=\Sigma^{i+j}$.  
   (The words ``up
    to weak equivalence'' are needed here because $i$~and~$j$ can be
    negative.) Given
    a map $f:X\to Y$ of asymmetric spectra, the fibre $F$ of $f$ is an
    asymmetric spectrum with $F_i$ given by the homotopy fibre of
    $X_i\to Y_i$; the cofibre $C$ is the asymmetric spectrum with
    $C_i$ given by the mapping cone of $X_i\to Y_i$.  It turns out that
    $\Sigma F$ and $C$ are naturally weakly equivalent, and that the
    cofibre of $Y\to C$ is naturally weakly equivalent to $\Sigma X$.
    This allows the homotopy category of asymmetric spectra to be
    given a triangulated structure, in which $F\to X\to Y$ or $X\to
    Y\to C$ are distinguished triangles. These distinguised triangles
    give long exact sequences on homotopy groups. 

    Any space gives rise
    to a suspension spectrum $\Sigma^\infty X$, where $(\Sigma^\infty
    X)_n=S^n\Smash X$; in some sense this is the free $S$-module
    on~$X$. The homotopy groups of $\Sigma^\infty X$ are the stable
    homotopy groups of~$X$.
  \end{numbered}

  In spite of these encouraging signs, the category
  of asymmetric spectra has one serious shortcoming: there is no
  obvious way to define an internal tensor product on the category.
  Given two asymmetric spectra $X$, $Y$, one would like to define
  $X\Tensor Y$ to be $X\aBoxtensor_SY$.
  There is a real difficulty in making a definition like this,
  stemming from the fact that $S$ is \emph{not} a commutative monoid
  with respect to $\aBoxtensor$. This is exactly the same difficulty
  that comes up in trying to form $M\Tensor_RN$ when $M$ and $N$ are
  left modules over the noncommutative ring~$R$, and the
  intention is that $M\Tensor_RN$ be another left $R$-module.
  The remedy for
  this is to provide a little extra structure in the underlying
  objects, enough structure so that $S$ becomes a commutative monoid.
\end{NumberedSubSection}

\begin{NumberedSubSection}{Symmetric spectra}
The trick is to build in symmetric group actions, hence the name,
\emph{symmetric spectra}.

\begin{defn}\label{SymmetricBoxtensor}
  A \emph{symmetric sequence} $X$ of  spaces is a collection
  $\{X_n\}_{n\ge0}$ of spaces, together with, for each $n$, a
  left action of the symmetric group $\Sigma_n$ on $X_n$. The (graded)
  product $X\Boxtensor Y$ of two such sequences is the sequence $Z$
  given by
  \[   Z_n = \coprod_{i+j=n} (\Sigma_n)_+
  \Smash_{\Sigma_i\times\Sigma_j} X_i\Smash Y_j\,.
  \]
\end{defn}
Here $(\Sigma_n)_+$ denotes the union of $\Sigma_n$ with a disjoint
basepoint. 

\begin{rem}\label{ElegantTensor}
  A symmetric sequence can equally well be thought of as a functor
  from the category of finite sets and isomorphisms to the category of
  spaces. From this point of view the graded product has the more
  elegant description
  \[
      (X\Boxtensor Y)(C) = \coprod_{A\cup B=C, \,\,A\cap B=\emptyset}
      X(A)\Smash Y(B)\,.
  \]
\end{rem}
  
The definition of $\Boxtensor$ differs from that of $\aBoxtensor$
(\ref{AsymmetricBoxtensor}) because
of the need to have symmetric group actions on the constituents of the
result. The product $\Boxtensor$ gives a symmetric monoidal structure
on the category of symmetric sequences of spaces; the unit for
$\Boxtensor$ is again the sequence $\epsilon$ mentioned after \ref{AsymmetricBoxtensor}, promoted to a
symmetric sequence in the only possible way.
The twist isomorphism $X\Boxtensor
Y\iso Y\Boxtensor X$ is composed, in the formulation from \ref{ElegantTensor}, of
the usual isomorphisms
\[
       X(A)\Smash Y(B) \iso Y(B)\Smash X(A)\,.
\]
In the formulation of \ref{SymmetricBoxtensor}, the twist isomorphism
combines isomorphisms $X_i\Smash Y_j\iso Y_j\Smash X_i$ with
right multiplication on $\Sigma_n$ by an element of the symmetric
group which conjugates
$\Sigma_i\times\Sigma_j$ to $\Sigma_j\times\Sigma_i$. 

The sequence $S$ described above after \ref{AsymmetricBoxtensor}
extends to a symmetric sequence in a natural way, where $\Sigma_n$
acts on $S^1\Smash\cdots\Smash S^1$ by permuting the factors. We will
denote this symmetric sequence by \Sphere, since it will correspond to
the sphere spectrum. The natural $\Sigma_i\times\Sigma_j$-equivariant
isomorphisms $S^i\Smash S^j\iso~S^{i+j}$ combine to give a natural map
$m:\Sphere\Boxtensor \Sphere\to \Sphere$. The key property of this map
is the following one.

\begin{lem}
The map $m$ gives $\Sphere$ the structure of a commutative monoid (with
respect to $\Boxtensor$) in the category of symmetric sequences of spaces.
\end{lem}

\begin{rem}
  Commutativity means that $m\tau=m$, where $\tau:\Sphere\Boxtensor \Sphere\to \Sphere\Boxtensor \Sphere$ is the twist isomorphism of the symmetric monoidal
  structure.
\end{rem}

The basic definitions are now clear.

\begin{defn}
  A \emph{symmetric spectrum} $X$ is a symmetric sequence of pointed
  spaces which is a left module over $\Sphere$
  (i.e., has been provided with a map $m_X: \Sphere\Boxtensor X\to X$ with
  appropriate unital and associativity properties).
\end{defn}

\begin{rem}
  Since $\Sphere$ is a commutative monoid, there is no real distinction
  between the notions of left and right modules; if $m_X:\Sphere\Boxtensor
  X\to X$ gives $X$ the structure of a left $S$-module, then
  $m_X\cdot\tau:X\Boxtensor \Sphere\to X$ gives $X$ the structure of a right
  $\Sphere$-module. Moreover, in this case $m_X\cdot\tau$ is a map of left
  $\Sphere$-modules. 
\end{rem}

\begin{defn}
  The \emph{tensor product} (or \emph{smash product}) of two symmetric
  spectra $X$ and $Y$ is the symmetric spectrum $X\Boxtensor_{\Sphere} Y$ defined
  by the coequalizer diagram
  \[
                X\Boxtensor \Sphere\Boxtensor Y \Rightarrow X\Boxtensor Y\to
                X\Boxtensor_{\Sphere} Y\,,
  \]
  where the two maps on the left are induced by the $\Sphere$-module
  structures of $X$ and $Y$. 
\end{defn}

From now on we will drop the word ``symmetric'' and call a symmetric
spectrum a spectrum.  The tensor product $X\Boxtensor_{\Sphere}Y$ is
denoted $X\Tensor Y$, or, in topological contexts, $X\Smash Y$. The
tensor product gives a symmetric monoidal structure on the category of
spectra; the unit for this structure is $\Sphere$.

\begin{defn}
  An \emph{\ring/} (or \emph{ring spectrum}) $R$ is a spectrum together with maps $S\to R$
  and $R\Tensor R\to R$ with appropriate unital and associativity
  properties. The \ring/ $R$ is commutative if the
  multiplication $m_R:R\Tensor R\to R$ is \emph{commutative}, i.e., if
  $m_R\cdot\tau=m_R$, where $\tau$ is the twist automorphism of
  $R\Tensor R$.
\end{defn}

We leave it to the reader to define modules over \aring/, tensor
products of modules, algebras over a commutative \ring, etc. Note that
\Sphere/ is \aring, that every spectrum is a
module over \Sphere, and that every \ring/ is, in fact, an algebra over
\Sphere. From this point of view \Sphere/ plays the role of the ground
ring for the category of spectra, just as \Z/ is the ground ring for
the category of chain complexes.

\begin{titled}{Mapping spectra}
Given two spectra $X$, $Y$, it is possible to define a mapping spectrum
$\Hom(X,Y)$. If the two spectra are modules over the \ring/ $R$,
there is also a spectrum $\HomR(X,Y)$ of $R$-module maps;
if in addition $R$ is commutative, $\HomR(X,Y)$  is an
$R$-module spectrum.
\end{titled}

\begin{titled}{Homotopy theory and derived constuctions}
  There is quite a bit of work to be done in setting up the homotopy
  theory of spectra; in particular, it is tricky to define the
  homotopy groups $\pi_iX$ of a spectrum $X$, or what comes to the
  same thing, to define weak equivalences (quasi--isomorphisms)
  between spectra \cite[\S3]{HSS}. Familiar issues of a homological
  algebra nature come up: for instance, tensor products or mapping
  spectra do not necessarily preserve weak equivalences unless the
  objects involved have freeness (cofibrancy) or injectivity
  (fibrancy) properties. These issues are handled in the non--additive
  context of spectra by Quillen's model category machinery
  \cite{rDSModel} \cite{rHoveyModel}, which essentially allows a great
  deal of homological algebra to be extended to sufficiently
  structured non--additive settings. Invoking this machinery leads to
  notions of derived tensor product and derived mapping spectrum.
\end{titled}

After the dust has settled, it is possible to prove that the homotopy
theory of spectra is equivalent to the homotopy theory of asymmetric
spectra. There is a shift operation as in \ref{AsymmetricShift}, as
well as a triangulated structure on the homotopy category.   This homotopy category is
obtained as usual from the category of spectra by inverting the weak
equivalences. Any
space~$X$ gives rise to a suspension (symmetric) spectrum, which we
continue to denote $\Sigma^\infty X$;  as in \ref{AsymmetricShift}, $(\Sigma^\infty X)_n$ is
$S^n\Smash X$, but now the symmetric group
acts on $S^n=S^1\Smash\cdots\Smash S^1$ by permuting the factors.

\begin{numbered}{Rings and DGAs vs. \rings/}
  The construction of \ref{AsymmetricChains} can be extended to
  convert any ring $R$ into \aring/ $R_{\Sphere}$ and any ordinary
  module over $R$ to a module over $R_{\Sphere}$. More generally, any
  chain complex over $R$ gives a module over $R_{\Sphere}$, and this
  correspondence provides an equivalence between the derived category
  of $R$ and the homotopy category of modules over $R_{\Sphere}$ (this
  is actually part of a Quillen equivalence between two model
  categories). In this paper we work consistently with \rings, and
  when an ordinary ring $R$ comes into play we do not usually
  distinguish in notation between $R$ and $R_{\Sphere}$. Note that
  $\pi_iR_{\Sphere}$ is $R$ if $i=0$, and $0$~otherwise.

  These considerations apply more generally if $R$ is a DGA
  \cite{ShipleyDGA}; there is an associated \ring/ $R_{\Sphere}$, and
  a Quillen equivalence which induces an equivalence between the
  derived category of $R$ (i.e. the category obtained from DG
  $R$-modules by inverting the quasi--isomorphisms) and the homotopy
  category of $R_{\Sphere}$-modules. Note that
  $\pi_iR_{\Sphere}=H_iR$.

  Any ring or DGA is an algebra over $\Z$, and the corresponding
  \ring/ $R_{\Sphere}$ is an algebra over $\Z_{\Sphere}$. The
  correspondence $R\mapsto R_{\Sphere}$ gives a bijection up to
  equivalence between DGAs and $\Z_{\Sphere}$-algebras, or between
  rings and $\Z_{\Sphere}$ algebras whose homotopy is concentrated in
  degree~$0$. (Actually, any \ring/ whose homotopy is concentrated in
  degree~0 is canonically a $\Z_{\Sphere}$-algebra, and so amounts to
  an ordinary ring.)

  The situation with commutativity is more complicated; commutative
  $\Z_{\Sphere}$-algebras correspond to $E_\infty$ DGAs (these are
  DGAs with a multiplication which is homotopy commutative up to
  explicit higher homotopies which are parametrized by the cells of an
  $E_\infty$~operad.) The prime example of such an $E_\infty$ algebra
  is the integral cochain algebra on a space $X$
  \cite[App.~C]{MandellPadic} (since all of our DGAs
  have differential which decreases dimension by~$1$, the cochain
  algebra is treated as a DGA by placing the $i$-dimensional cochains
  in degree~$-i$). From a homotopical point of view, $C^*(X;\Z)$ is
  $\Hom_{\Sphere}(\Sigma^\infty X_+, \Z_{\Sphere})$, where $X_+$ is $X$ with a
  disjoint basepoint adjoined. The commutative \ring/ structure is
  derived from the multiplication on $\Z_{\Sphere}$ and the diagonal
  map on~$X$.

  To repeat, commutative $\Z_{\Sphere}$-algebras do
  \emph{not} correspond to commutative DGAs. However, there is one
  bright note: commutative $Z_{\Sphere}$-algebras with homotopy
  concentrated in degree~$0$ do correspond bijectively up to
  equivalence to commutative rings.
\end{numbered}

\begin{titled}{Terminological caveat}
  In the literature, \rings/ are sometimes styled \emph{$A_\infty$
    ring spectra} (or \emph{structured ring spectra}) and commutative
  \rings/ \emph{$E_\infty$ ring spectra}. The terms \emph{ring
    spectrum} and \emph{commutative ring spectrum} are occasionally
  used even today for a much weaker notion in which various diagrams
  involving the multiplication map are only required to commute up to
  homotopy.
\end{titled}
\end{NumberedSubSection}

\section{Some basic constructions with modules}\label{CModuleProperties}

This section looks into some constructions with \rings/ and modules
which we refer to in the rest of the paper. We first describe some
Postnikov constructions which allow modules to be filtered in such a
way that the successive (co)fibres are ``Eilenberg-MacLane objects'',
in the sense that they have nonvanishing homotopy in only a single
dimension (\ref{ConnectivePostnikov}, \ref{CoconnectivePostnikov}).
Next, we show that in many cases these Eilenberg-MacLane objects are
determined by the homotopy group which appears
(\ref{UniqueModuleStructures}), although there are surprises
(\ref{ExoticEilenbergMacLane}). We end by formulating
``directionality'' properties for modules (\ref{UpAndDown}), and
studying finiteness conditions (\ref{FinitelyBuilt}).

First, some terminology.
\Aring/ $R$ is \emph{connective} if $\Hof iR=0$ for
      $i<0$ and \emph{coconnective} if $\Hof iR=0$ for $i>0$. An
      $R$-module $M$ is \emph{bounded below} if $\Hof iM=0$ for
      $i<\!\!<0$, and \emph{bounded above} if $\Hof iM=0$ for
      $i>\!\!>0$.

\begin{NumberedSubSection}{Postnikov constructions}
There are subtle differences between the connective and coconnective
cases. 

\begin{prop}\label{ConnectivePostnikov}
  Suppose that $R$ is connective, that $M$ is an $R$-module, and that
  $n$ is an integer. Then there is a natural $R$-module $P_nM$ with
  $\Hof i(P_nM)\iso 0$ for $i>n$,
  together with a natural map $M\to P_nM$ inducing isomorphisms on
  $\Hof i$ for $i\le n$.
\end{prop}

\begin{prop}\label{CoconnectivePostnikov}
  Suppose that $R$ is coconnective with $\Hof0R$ a field, that $M$ is
  an $R$-module, and that $n$ is an integer. Then there is an
  $R$-module $Q_nM$ with $\Hof i(Q_nM)=0$ for $i<n$, together with
  a map $M\to Q_nM$ inducing isomorphisms on $\Hof i$ for
  $i\ge n$.
\end{prop}

\begin{rem}
  In the above situations there are maps $P_nM\to P_{n-1}M$ or
  $Q_{n}M\to Q_{n+1}M$ inducing isomorphisms on all appropriate nonzero
  homotopy groups. In the first case $M\weq\holim_n P_nM$, while in
  the second $M\weq \holim_nQ_nM$. The fibre of the map
  $P_nM\to P_{n-1}M$ or $Q_{n}M\to Q_{n+1}M$ is a module with only one
  nonvanishing homotopy group.
\end{rem}

\begin{rem}
  The construction of $Q_nM$ cannot be made functorial in any
  reasonable sense. Consider the DGA $\EE$ of
  \cite[\S3]{dwyer-greenlees}; $\EE$ is coconnective and
  $\Hof0{\EE}\iso\Fp$. Then $\Hof0\HomE(\EE,\EE)\iso\Hof0\EE\iso\Fp$,
  while $\Hof0\HomE(Q_0\EE,Q_0\EE)=\pi_0\HomE(\Fp,\Fp)\weq\Zpadic$. We
  are using topological notation: $\Zpadic$ is the ring of $p$-adic
  integers. Since there is no
  additive map $\Fp\to\Zpadic$, there is no way to form $Q_0\EE$
  functorially from $\EE$.
\end{rem}

\begin{rem}\label{DefineHoCoLim} We have used the notion of ``$\holim$'' above, and later
  on we will use ``$\hocolim$''. If $M_0\to M_1\to \cdots$ is an
  inductive system of spectra, then $\hocolim_nM_n$ is the fibre of
  $1-\sigma$, where $\sigma:\coprod M_n\to\coprod M_n$ is the shift
  map. There are isomorphisms
  $\pi_i\hocolim_nM_n\iso\colim_n\pi_iM_n$.  Dually, if $\cdots \to
  N_1\to N_0$ is an inverse system of spectra, $\holim_nN_n$ is
  defined to be the fibre of $1-\sigma$, where $\sigma:\prod
  N_n\to\prod N_n$ is the shift map. In this case there are short
  exact sequences
  \[
        0 \to \lim_n{}\!^1\,\pi_{i+1}N_n\to\pi_i\holim_nN_n\to \lim_n \pi_iN_n\to
        0\,.
  \]
  Homotopy colimits and limits are defined for arbitrary small
  diagrams of spectra \cite[18.1]{rHirschhorn}, but we will not need
  them in this generality.
\end{rem}

For the proofs we need to make attaching constructions.

  \begin{defn}\label{DefineAttach}
    Suppose that $R$ is \aring/ and that $X$ and $Y$ are $R$-modules.
    Then $Y$ is obtained from $X$ by
    \emph{attaching an $R$-module $A$} if there is a
    cofibration sequence $A\to X\to Y$. If $\{A_\alpha\}$ is a
    collection of $R$-modules, then $Y$ is obtained from $X$ by
    \emph{attaching copies of the $A_\alpha$} if there is a cofibration
      sequence $U\to X\to Y$ in which $U$ is a coproduct of elements
      from $\{A_\alpha\}$.
    Even more
    generally, $Y$ is obtained from $X$ by \emph{iteratively attaching
      copies of $\{A_\alpha\}$} if $Y$ is the colimit of a
    directed system $Y_0\to Y_1\to\cdots$, such that
    $X_0=X$ and
    $X_{n+1}$ is obtained from
    $X_n$ by attaching copies of the $A_\alpha$.
  \end{defn}

\begin{titled}{Proof of \ref{ConnectivePostnikov}}
  Form $P_nM$ by iteratively attaching copies of $\Shift^i R$, $i>n$
  to $M$ (\ref{DefineAttach}) to kill off the homotopy of $M$ above
  dimension~$n$. More specifically, for each element $x\in\pi_{i_x}M$
  with $i_x>n$, construct a map $\Sigma^{i_x}R\to M$ which takes the
  unit in $\pi_0R$ to~$x$. Let $C(M)$ be the cofibre of the resulting map
  $\coprod_x\Sigma^{i_x}R\to M$, and observe that the map $M\to C(M)$
  induces isomorphisms on $\pi_i$ for $i\le n$ and is zero on $\pi_i$
  for $i>n$. Repeat the process, and let $P_nM=\hocolim_kC^k(M)$,
  where $C^k(M)=C(C^{k-1}M)$.
  The construction can be made functorial by doing the
  attachments over all maps with domain $\Sigma^iR$, $i>n$, and not
  bothering to choose representative maps from each homotopy class. \qed
\end{titled}

\begin{titled}{Proof of \ref{CoconnectivePostnikov}}
  Given an $R$-module $X$ and an integer $m$, choose a basis for $\Hof mX$ 
  over $\Hof0R$, and let $V_mX$ be a sum of copies of $\Shift ^mR$,
  one for each basis element. There is a map $V_mX\to X$ which
  induces an isomorphism on $\Hof m$. Let $C(M)$ be the cofibre of the
  map $\coprod_{m<n}V_mX\to X$, and observe that the map $M\to C(M)$
  induces an isomorphism on $\pi_i$ for $i\ge n$, and is zero on
  $\pi_i$ for $i<n$. Repeat the process, and let $Q_nM=\hocolim_k C^k(M)$.

  It is the fact that $\Hof0R$ is a field which guarantees that the
  attachment producing $C(M)$ can be done without introducing new
  homotopy in dimensions $\ge n$. However, the attachment must be done
  minimally, and it is this requirement that prevents the construction
  from being functorial. \qed
\end{titled}

\end{NumberedSubSection}

\begin{NumberedSubSection}{Uniqueness of module structures}
We first aim for the following elementary uniqueness result.

\begin{prop}\label{UniqueModuleStructures}
  Suppose that $R$ is connective or that $R$ is coconnective with
  $\Hof0R$ a field, and that $M$ and $N$ are $R$-modules with
  nonvanishing homotopy only in a single dimension $n$. Then $M$ and
  $N$ are equivalent as $R$-modules if and only if $\Hof nM$ and 
  $\Hof nN$ are isomorphic over $\Hof 0R$.
\end{prop}

\begin{rem}\label{EMObjects}
  It follows easily from the proof below that if $R$ is as in
  \ref{UniqueModuleStructures}, $A$ is a discrete module over $\Hof0R$, and $n$
  is an integer, then there exists an
  $R$-module $K(A,n)$ with $\Hof nK(A,n)\iso A$ (over $\Hof 0R$) and
  $\Hof iK(A,n)\iso 0$ for $i\ne n$. If $R$ is connective the
  construction of $K(A,n)$ can be made functorially in $A$, otherwise
  in general not. If $A$ and $B$ are two discrete $\Hof0R$-modules, the natural
  map \[\Hof0\HomR(K(A,n),K(B,n))\to\hom_{\Hof0R}(A,B)\] is an
  isomorphism if $R$ is connective but only a surjection in general if
  $R$ is coconnective and $\pi_0R$ is a field.
\end{rem}

\begin{rem}\label{ExoticEilenbergMacLane}
  A startling cautionary note is struck by the fact that
  if $R$ is coconnective and $\Hof0R$ is not a field,
  the conclusion of \ref{UniqueModuleStructures} is \emph{not}
  necessarily true. We sketch an example of \aring/ $R$ and two
  $R$-modules $M$ and $N$ with homotopy concentrated in degree~$0$,
  such that $\pi_0M\iso\pi_0N$ as modules over $\pi_0R$, but $M$ is
  not equivalent to $N$ as an $R$-module. Let $S$ be the ring $\Z[t]$;
  make $\Z$ into a discrete $S$-module by letting $t$ act as
  multiplication by zero, and $\Z/p^\infty$ into a discrete $S$-module
  by letting $t$ act by multiplication by~$p$. Let $\Fp[t]/t^\infty$
  be the discrete $S$-module $\Fp[t,t^{-1}]/\Fp[t]$.  We let
  $R=\End_S(\Z)$, $M=\Hom_S(\Z,\Z/p^\infty)$, and
$N=\Hom_S(\Z,\Fp[t,t^{-1}]/\Fp[t])$.  Then $\pi_*R=\Ext^{-*}_S(\Z,\Z)$
is an exterior algebra over~$\Z$ on a generator of degree~$-1$; $\pi_*M =
\Ext^{-*}_S(\Z,\Z/p^\infty)$ is a copy of $\Z/p$ in degree~$0$, and
$\pi_*N$ is isomorphic to~$\pi_*M$.  Both $M$ and $N$ are right
$R$-modules, but they are not equivalent as right $R$-modules, since for instance by
\cite[2.1]{dwyer-greenlees} (cf. \ref{TMoritaTheorem}) there are
equivalences
  \[
      M\Tensor_R \Z \weq \Z/p^\infty\text{\quad and\quad} N\Tensor_R\Z\weq\Fp[t]/t^\infty\,.\]
\end{rem}

\begin{titled}{Proof of \ref{UniqueModuleStructures}}
  One way to prove this is to construct a suitable spectral sequence
  converging to $\Hof*\HomR(M,N)$; under the connectivity assumptions
  on $R$, $\hom_{\Hof 0R}(\Hof nM, \Hof nN)$ will appear in one corner
  of the  $E_2$-page and subsequently remain undisturbed 
  for positional reasons. This implies that any map $\Hof nM\to
  \Hof nN$ of $\Hof0R$-modules, in particular any isomorphism, can
  be realized by an $R$-map $M\to N$. We will take a more elementary
  approach. Assume without loss of generality that $n=0$ and suppose
  that there are isomorphisms $\Hof 0M\iso \Hof 0N\iso A$ over $\Hof 0R$. First we treat the
  case in which $R$ is connective. Find a free presentation
  \[ \phi_1\to \phi_0\to A\to 0\]
  of $A$ over $\Hof0R$ and construct a map $F_1\to F_0$ of $R$-modules
  such that each $F_i$ is a sum of copies of $R$, and such that $\Hof0
  F_1\to \Hof 0F_0$ is $\phi_1\to\phi_0$. Let $C$ be the cofibre
  of $F_1\to F_0$. By inspection $\Hof 0C\iso A$ and there are
  isomorphisms $\Hof 0\HomR(C,M)\iso\hom_{\Hof 0R}(A,A)$ and
  $\Hof0\HomR(C,N)\iso\hom_{\Hof 0R}(A,A)$. Choose maps $C\to M$ and
  $C\to N$ which induce isomorphisms on $\Hof0$, and apply the functor
  $P_0$ (\ref{ConnectivePostnikov}) to conclude $M\weq N$. Now suppose
  that $R$ is coconnective, and that $\Hof0R$ is a field. Write $A\iso
  \oplus_\alpha \Hof0R$ over $\Hof0R$, let $F=\oplus_\alpha R$, and
  construct maps $F\to M$ and $F\to N$ inducing isomorphisms on $\Hof
  0$. Consider $Q_0F$ (\ref{CoconnectivePostnikov}). Since $Q_0F$ is
  obtained from $F$ by attaching copies of $\Shift ^iR$, $i<0$, there
  are surjections (not necessarily isomorphisms) $\Hof
  0\HomR(Q_0F,M)\twoheadrightarrow\hom_{\Hof 0R}(A,A)$ and
  $\Hof0\HomR(Q_0F,N)\twoheadrightarrow\hom_{\Hof 0R}(A,A)$.  Clearly,
  then, there are equivalences $Q_0F\to M$ and $Q_0F\to N$. \qed
\end{titled}

\end{NumberedSubSection}

\begin{NumberedSubSection}{Upward and downward (finite) type}\label{UpAndDown}
  Suppose that $M$ is an $R$-module. We say that $M$ is \emph{of
    upward type} if there is some integer $n$ such that up to
  equivalence $M$ can be built by starting with the zero module and
  iteratively attaching copies of $\Shift^i R$, $i\ge n$; $M$ is of \emph{upward
    finite type} if the construction can be done in such a way that
  for any single $i$ only a finite number of copies of $\Shift ^iR$ are
  employed. Similarly, $M$ is \emph{of downward type} if there is some
  integer $n$ such that $M$ can be built by starting with the zero
  module and iteratively attaching copies of $\Shift^iR$, $i\le n$; $M$ is of
  \emph{downward finite type} if the construction can be done in such
  a way that for any single $i$ only a finite number of copies of
  $\Shift^iR$ are employed.

  We look for conditions under which an $R$-module has upward or
  downward (finite) type.

\begin{prop}\label{UpwardFiniteType}
  Suppose that $R$ is a connective \ring, and that $M$ is a module
  over $R$ which is bounded below.  Then $M$ is of upward type.  If in
  addition $\Hof 0R$ is Noetherian and the groups $\Hof iR$ and $\Hof iM$ ($i\in\Z$) are individually finitely generated over $\Hof0R$,
  then $M$ is of upward finite type.
\end{prop}

\begin{proof}
  Suppose for definiteness that $\pi_iM=0$ for $i<0$. We inductively
  construct maps $X_n\to M$ such that $\pi_iX_n\to\pi_iM$ is an
  isomorphism for $i<n$ and an epimorphism for $i=n$, and such that
  $X_n$ is constructed from~$0$ by attaching copies of
  $\Sigma^iR$ for $0\le i\le n$. Let $X_{-1}=0$. Given $X_n\to M$,
  construct a surjection
  \[
      \coprod_\alpha \pi_0R\twoheadrightarrow\ker(\pi_nX_n\to\pi_nM)
  \]
  of modules over $\pi_0R$, realize this surjection by a map
  $\coprod_\alpha\Sigma^nR\to X_n$, and let $X_n'$ be the cofibre of
  this map. The map $X_n\to M$ extends to a map $X_n'\to M$ which is
  an isomorphism on $\pi_i$ for $i\le n$. Choose a surjection
  $\coprod_\beta\pi_0R\to \pi_{n+1}M$ of $\pi_0R$-modules, realize this by
  a map $Y=\coprod_\beta\Sigma^{n+1}R\to M$, and let
  $X_{n+1}=X_n'\coprod Y$. The map $X_n'\to
  M$ then extends to a map $X_{n+1}\to M$ with the desired properties.
  A homotopy group calculation gives that $\hocolim_n X_n\weq M$, and
  so it is clear that $M$ is of upward type over $R$. Under the stated
  finiteness assumptions, it is easy to prove inductively that the
  homotopy groups of $X_n$ are finitely generated over $\pi_0R$, and
  consequently that
  the above coproducts of suspensions of $R$ 
  can be chosen to be finite.
\end{proof}

\begin{prop}\label{DownwardFiniteType}
  Suppose that $R$ is a coconnective \ring/ such that $\Hof0R$ is a
  field, and that $M$ is an $R$-module which is bounded above. Then
  $M$ is of downward type. If in addition $\Hof {-1}R=0$ and the
  groups $\Hof iR$ and $\Hof iM$ ($i\in\Z$) are individually finitely
  generated over $\Hof0R$, then $M$ is of downward finite type.
\end{prop}

\begin{proof}
  Given an $R$-module $X$, let $V_mX$ be as in the proof of
  \ref{CoconnectivePostnikov}
  and let $W_mX$ be the cofibre
  of $V_mX\to X$. Now suppose that $M$ is nontrivial and bounded above,
  let $n$ be the greatest integer such that $\Hof nM\ne0$, and let
  $W_nM$ be the cofibre of the map $V_nM\to M$. Iteration 
  gives a sequence of maps $M\to W_nM\to W_n^2M\to\cdots$, and we let
  $W^\infty_nM=\hocolim_k W^k_nM$. Then $\Hof nW^\infty_nM=\colim_k\Hof nW_n^kM=0$.
  Define modules $U_i$ inductively by $U_0=M$, $U_{i+1}
  = W^\infty_{n-i}U_i$. There are maps $U_i\to U_{i+1}$ and it is
  clear that $\hocolim U_i\weq0$. Let $F_i$ be the homotopy fibre of
  $M\to U_i$. Then $\hocolim F_i$ is equivalent to $M$, and $F_{i+1}$
  is obtained from $F_{i}$ by repeatedly attaching copies of
  $\Shift^{n-i-1}R$.  This shows that $M$ is of downward type. If
  $\Hof{-1}R=0$, then $\Hof {n-i}W_{n-i}U_i\iso0$, so that $W_{n-i}^\infty U_i\weq W_{n-i}U_i$.
  Under the stated finiteness hypotheses, one sees by an inductive
  argument that the groups $\Hof jU_i$, $j\in\Z$, are finite
  dimensional over $k$, and so $F_{i+1}$ is obtained from $F_i$ by
  attaching a finite number of copies of $\Shift^{n-i-1}R$.  This
  shows that $M$ is of downward finite type.
\end{proof}

\end{NumberedSubSection}

\begin{NumberedSubSection}{(Finitely) built}\label{FinitelyBuilt}
  A subcategory of \RMod/ is \emph{thick} if it is 
  closed under equivalences, triangles, and retracts; here closure
  under triangles means that given any distinguished
  triangle (i.e. cofibration sequence) with two of its terms in the
  category, the third is in the category as well. The subcategory is
  \emph{localizing} if in addition it is closed under arbitrary
  coproducts (or equivalently, under arbitrary homotopy colimits).
  If $A$ and $B$ are $R$-modules, we say that $A$ is
  \emph{finitely built from~$B$} if $A$ is in the smallest thick
  subcategory of \RMod/ which contains $B$; $A$ is
  \emph{built from~$B$} if it is contained in the smallest localizing
  subcategory of \RMod/ which contains $B$.

  Given an augmented $k$-algebra $R$, we look at the question of when
  an $R$-module $M$ is (finitely) built from $R$ itself or from $k$.
  We have already touched on related issues. It is in fact not hard to see
  that \emph{any} $R$-module is built from $R$; propositions
  \ref{UpwardFiniteType} and \ref{DownwardFiniteType} amount to
  statements that sometimes this building can be done in a controlled
  way.

\begin{prop}\label{BestOfAllCoregularity}
  Suppose that $k$ is a field, that $R$ is an augmented $k$-algebra,
  and that $M$ is an $R$-module.  Assume either that $R$ is connective
  and the kernel of the augmentation $\Hof 0R\to k$ is a nilpotent
  ideal , or that $R$ is coconnective and $\Hof 0R\iso k$. Then an
  $R$-module $M$ is finitely built from $k$ over $R$ if and only if
  $\Hof *M$ is finite dimensional over $k$.
\end{prop}

\begin{rem}\label{BuiltFromk}
  A similar argument shows that if $R$ is coconnective and $\Hof0R\iso k$,
  then any $R$-module $M$ which is bounded below is built from $k$
  over $R$. It is only necessary to note that the fibre $F_n$ of $M\to
  Q_nM$ is built from $k$ (it has only a finite number of nontrivial
  homotopy groups) and that $M\weq \hocolim F_n$. Along the same lines,
  if $R$ is connective and $\Hof0R$ is as in
  \ref{BestOfAllCoregularity}, then any $R$-module $M$ which is
  bounded above is built from $k$ over $R$.
\end{rem}

\begin{titled}{Proof of \ref{BestOfAllCoregularity}}
  It is clear that if $M$ is finitely built from $k$ then $\Hof*M$ is
  finite dimensional. Suppose then that $\Hof*M$ is
  finite-dimensional, so that in particular $\Hof iM$ vanishes for all
  but a finite number of $i$. By using the Postnikov constructions $P_*$
  (\ref{ConnectivePostnikov}) or $Q_*$ (\ref{CoconnectivePostnikov}),
  we can find a finite filtration of $M$ such that the successive
  cofibres are of the form $K(\Hof nM,n)$ (\ref{EMObjects}). It
  is enough to show that if $A$ is a discrete module over $\Hof0R$ which is
  finite-dimensional over $k$, then $K(A,n)$ is finitely built from
  $k$ over $R$. But this follows from \ref{EMObjects} and that fact
  that under the given assumptions, $A$ has a finite filtration by
  $\Hof0R$-submodules such that the successive quotients are
  isomorphic to $k$. \qed
\end{titled}

We need one final result, in which $k$ and $R$ play reciprocal roles.
In the following proposition, there is a certain arbitrariness to the choice of which \ring/ 
is named \EE/ and which is named~$R$; we have
picked the notation so that the formulation is parallel to
\ref{DownwardFiniteType}.

\begin{prop}\label{SpecialDownwardFiniteType}\label{DownwardSmall}
  Suppose that $k$ is a field, and that $\EE$ is a connective
  augmented $k$-algebra such that the kernel of the augmentation
  $\pi_0\EE\to k$ is a nilpotent ideal. Let $M$ be an \EE-module, let
  $N=\Hom_{\EE}(M,k)$, and let $R=\End_ {\EE}(k)$. If $M$ is finitely
built from $k$ over \EE/ (i.e. $\pi_*M$ is finite-dimensional
over~$k$) then $N$ is finitely built from $R$ over $R$.  If $M$ is
bounded below and each $\pi_iM$ is finite-dimensional over~$k$, then
$N$ is of downward finite type over~$R$.
\end{prop}

\begin{rem} This proposition
  will let us derive the conclusion of \ref{DownwardFiniteType} in
  some cases in which $\pi_{-1}R\ne0$.
\end{rem}

\begin{titled}{Proof of \ref{SpecialDownwardFiniteType}}
  Suppose that $X$ is some $\EE$-module.
  It is elementary that if $M$ is finitely built from~$X$ over \EE,
  then $\HomE(M,k)$ is finitely built from $\HomE(X,k)$ over
  $\End_{\EE}(k)$. Taking $X=k$ gives the first statement of the
proposition.  

Suppose then that $M$ is bounded below, and that each $\pi_iM$ is
finite-dimensional over~$k$. Let $M_i$ denote the Postnikov stage
$P_iM$ (\ref{ConnectivePostnikov}), so that $M$ is equivalent to
$\holim_iM_i$, and let $N_i=\HomE(M_i,k)$. We claim that $N$ is equivalent to
$\hocolim_iN_i$.  This follows from the fact that $M_n$ is
obtained from $M$ by attaching copies of $\Shift^i\EE$ for $i>n$, and
so the natural map $N_n\to N$ induces isomorphisms
on $\Hof i$ for $i\ge-n$. 
  The triangle $M_n\to M_{n-1}\to K(\Hof nM, n+1)$
(\ref{EMObjects}) dualizes to give a triangle 
\[\HomE(K(\Hof nM,n+1),
k)\to N_{n-1}\to N_n\,.\] Since (as in the proof of
\ref{BestOfAllCoregularity}) $\pi_nM$ has a finite filtration by
$\pi_0\EE$-submodules in which the sucessive quotients are isomorphic
to~$k$, it follows that that $N_n$ is obtained from $N_{n-1}$ by
attaching a finite number of copies of $\HomE(\Shift^{n+1}k,
k)\weq\Shift^{-(n+1)}R$.  Since $M_i\weq 0$ for $i<\!\!<0$, the
proposition follows.  \qed
\end{titled}

\end{NumberedSubSection}

\section{Smallness}\label{CContext}

In this section we describe the main setting that we work in.
We start with a pair $\Pair Rk$, where $R$ is \aring/ and $k$ is
an $R$-module.  Eventually we assume that $k$ is an
$R$-module via \aring/ homomorphism $R\to k$.

We begin by discussing cellularity (\ref{DiscussCellularModules}) and
then describing some smallness hypotheses under which cellular
approximations are given by a simple formula (\ref{TLatentMorita}).
These smallness hypotheses lead to various homotopical formulations of
smallness for \aring/ homomorphism $R\to k$ (\ref{DefineRegular}).
Let $\EE=\EndR(k)$. We show that the smallness conditions have a
certain symmetry under the interchange $R\leftrightarrow\EE$, at least
if $R$ is complete in an appropriate sense (\ref{RegularCriterion}),
and that the smallness conditions also behave well with respect to
``short exact sequences'' of \rings/
(\ref{EhregularInTheMiddle}). Finally, we point out that in some
algebraic situations the notion of completeness from
\ref{RegularCriterion} amounts to ordinary completeness with respect
to powers of an ideal, and that in topological situations it
amounts to convergence of the Eilenberg-Moore spectral sequence (\ref{ChainCochainStuff}).

\begin{NumberedSubSection}{Cellular modules}\label{DiscussCellularModules} 
  A map $U\to V$ of $R$-modules is a \emph{$k$-equivalence} if the
  induced map $\HomR(k,U)\to\HomR(k,V)$ is an equivalence.  An
  $R$-module $M$ is said to be \emph{$k$-cellular} or
  \emph{$k$-torsion} (\cite[\S4]{dwyer-greenlees}, \cite{Farjoun}) if
  any $k$-equivalence $U\to V$ induces an equivalence
  $\HomR(M,U)\to\HomR(M,V)$.   A $k$-equivalence between $k$-cellular objects is necessarily an
  equivalence. 

  It it not hard to see that any $R$-module which is built from $k$ in
  the sense of \ref{FinitelyBuilt} is $k$-cellular, and in fact it
  turns out that an $R$-module is $k$-cellular
  if and only if
  it is built from $k$ (cf.
  \cite[5.1.5]{rHirschhorn}). The proof involves using a version
  of Quillen's small object argument to show that for any $R$-module
  $M$, there exists a
  $k$-equivalence $M'\to M$ in which $M'$ is built from~$k$. If $M$ is $k$-cellular, this
  $k$-equivalence must be an equivalence.
  
We let 
$\DCellOfPair Rk$ denote the full subcategory
  of the derived category $\Derived(R)$.
  containing the $k$-cellular objects.
  For any $R$-module $X$ there is a $k$-cellular object $\CellOf k(X)$
  together with a $k$-equivalence $\CellOf k(X)\to X$; such an object
  is unique up to a canonical equivalence and is called the
  \emph{$k$-cellular approximation} to $X$. If we want to emphasize
  the role of $R$ we write $\CellOverOf Rk(X)$.

  \begin{rem}
    If $R$ is a commutative ring and $k=R/I$ for a finitely generated
    ideal~$I\subset R$, then an $R$-module $X$ is $k$-cellular if and
    only if each element of $\pi_*X$ is annihilated by some power
    of~$I$ \cite[6.12]{dwyer-greenlees}. The chain complex incarnation
    (\ref{OurNotation}) of $\CellOverOf Rk(X)$ is the local
    cohomology object $\mathbf{R}\Gamma_I(R)$
    \cite[6.11]{dwyer-greenlees}.
  \end{rem}

We are interested in a particular approach to constructing
$k$-cellular appproximations, and it is convenient to have some
terminology to describe it.

\begin{defn}\label{EffectivelyConstructible}
  Suppose that $k$ is an $R$-module and that $\EE=\EndR(k)$.
  An $R$-module  $M$ is said to be \emph{effectively
    constructible from $k$} if the natural map
  \[
       \HomR(k,M)\Tensor_{\EE}k\to M
   \]
  is an equivalence.
\end{defn}

Note that $\HomR(k,M)\Tensor_{\EE}k$ is always $k$-cellular over~$R$,
because $\HomR(k,M)$ is $\EE$-cellular over~$\EE$
(\ref{FinitelyBuilt}). The following
lemma is easy to deduce from the fact that the map $\Cell_k(M)\to M$
is a $k$-equivalence.

\begin{lem}\label{ThreeConstructible}
  In the situation of \ref{EffectivelyConstructible}, the following
  conditions are equivalent:
  \begin{enumerate}
  \item $\HomR(k,M)\Tensor_{\EE}k\to M$ is a $k$-equivalence.
  \item $\HomR(k,M)\Tensor_{\EE}k\to M$ is a $k$-cellular
    approximation.
   \item $\Cell_k(M)$ is effectively constructible from $k$.
  \end{enumerate}
\end{lem}

\end{NumberedSubSection}

\begin{NumberedSubSection}{Smallness}
In the context above of $R$ and~$k$, we will consider three finiteness
conditions derived from the notion of ``being finitely built'' (\ref{FinitelyBuilt}).

\begin{defn}\label{DefineSmallness}
  The $R$-module $k$ is \emph{\smalll/} if $k$ is finitely built from
  $R$, and \emph{\cosmall/} if $R$ is finitely built from $k$.
  Finally, $k$ is \emph{\mixedsmall/} if there exists an $R$-module
  $K$, such that $K$ is finitely built from $R$, $K$ is finitely built
  from $k$, 
  and $K$ builds~$k$.
   The object $K$ is then called
  a \emph{Koszul complex} associated to $k$ (cf. \ref{OrdinaryRingERegular}).
\end{defn}

\begin{rem}
The $R$-module $k$ is \smalll/ if and only if $\HomR(k,\whatever)$
commutes with arbitrary coproducts; if $R$ is a ring this is equivalent
to requiring that $k$ be a perfect complex, i.e., isomorphic in
$\Derived(R)$ to a chain complex of finite length whose constituents
are finitely generated projective $R$-modules.
\end{rem}

\begin{rem}\label{LatentMeansSameCells}
  The condition in 
  \ref{DefineSmallness} that $k$ and $K$ can be built from one another
  implies that an $R$-module $M$ is built from $k$ if and only if it
  is built from $K$; in particular, $\DCellOfPair Rk=\DCellOfPair RK$.
   If $k$ is either \smalll/ or \cosmall/ it is also
  \mixedsmall; in the former case take $K=k$ and in the latter $K=R$.

  In \cite{dgi:algebra} we explore the relevance of the concept of
  proxy--smallness to commutative rings.
\end{rem}

One of the main results of
\cite{dwyer-greenlees} is the following; although in
\cite{dwyer-greenlees} it is phrased for \DGAs, the proof for general
\rings/ is the same. If $k$ is an $R$-module, let $\EE=\EndR(k)$,
let $\E$ be the functor
which assigns to an $R$-module $M$ the right $\EE$-module
$\HomR(k,M)$, and let $\T$ be the functor which assigns to a right
$\EE$-module $X$ the $R$-module $X\Tensor_{\EE}k$. 

\begin{thm}\label{TMoritaTheorem}\label{TMoritaTheoremComments}\label{MoritaInverses}
  \cite[2.1, 4.3]{dwyer-greenlees} If $k$ is \asmalll/ $R$-module, then the functors $\E$ and $\T$
  above induce adjoint categorical equivalences
  \[
      \T : \Derived(\EE\op) \leftrightarrow \DCellOfPair Rk :\EE\,. 
  \]
  All $k$-cellular $R$-modules are effectively constructible from~$k$.
\end{thm}

There is a partial generalization of this to the \mixedsmall/ case.

\def\1{J}
\begin{thm}\label{TLatentMorita}\label{TLatentMoritaEquiv}
  Suppose that $k$ is \amixedsmall/ $R$-module with Koszul complex
  $K$.  Let $\EE=\EndR(k)$, $\1 =\HomR(k,K)$, and $\EEK=\EndR(K)$.
  Then the three categories
  \begin{center}
   $\DCellOfPair Rk$, 
   $\DCellOfPair {\EE\op}\1$, 
   $\Derived(\EEK\op)$
  \end{center} are all equivalent to one another.
  All $k$-cellular $R$-modules are effectively constructible from~$k$.
\end{thm}

\begin{rem}
  We leave it to the reader to inspect the proof below and write down the functors that induce the
  various categorical equivalences.
\end{rem}

\begin{titled}{Proof of \ref{TLatentMorita}}
  We will show that $\1$ is \asmalll/ $\EE\op$-module, and that the
  natural map $\EEK\to\EndOver{\EE\op}(\1 )$ is an equivalence. The
  first statement is then proved by applying \ref{TMoritaTheorem} serially to
  the pairs $\Pair{\EEK\op}{\1}$ and $\Pair RK$ whilst keeping
  \ref{LatentMeansSameCells} in mind.  For the smallness, observe that
  since $K$ is finitely built from $k$ as an $R$-module, $\1
  =\HomR(k,K)$ is finitely built from $\EE=\HomR(k,k)$ as a right
  $\EE$-module. Next, consider all $R$-modules $X$ with the property that
  for any $R$-module $M$ the natural map
  \begin{equation}\label{RecoverMaps}
          \HomR(X,M)\to \Hom_{\EE\op}(\HomR(k,X), \HomR(k,M ))
  \end{equation}
  is an equivalence. The class includes $X=k$ by inspection, and hence
  by triangle arguments any $X$ finitely built from $k$, in
  particular $X=K$.   

  For the second statement, suppose that $M$ is $k$-cellular. By
  \ref{LatentMeansSameCells}, $M$ is also $K$-cellular and hence
  (\ref{MoritaInverses}) effectively constructible from~$K$. In other
  words, the natural map
  \[\HomR(K,M)\TensorOver{\EEK}K\to M\] is an equivalence.
  We wish to analyze the domain of the map. As above (\ref{RecoverMaps}),
  $\HomR(K,M)$ is equivalent to
  $\Hom_{\EE\op}(\1,\HomR(k,M))$, which, because $\1$ is \smalll/ as a
  right $\EE$-module, is itself equivalent to
  $\HomR(k,M)\TensorOver{\EE}\Hom_{\EE\op}(J,\EE)$. Since
  $\EE\weq\HomR(k,k)$, the second factor of the tensor product is
  (again as with \ref{RecoverMaps}) equivalent to $\HomR(K,k)$. We
  conclude that the natural map
  \[
      \HomR(k,M)\TensorOver{\EE}(\HomR(K,k)\TensorOver{\EEK}K)\to M
  \]
  is an equivalence. But the factor
  $\HomR(K,k)\TensorOver{\EEK}K$ is equivalent to $k$, since by
  \ref{TMoritaTheorem} the $K$-cellular module~$k$ is effectively
  constructible from~$K$. Hence $M$ is effectively constructible from~$k$.
\qed
\end{titled}

\end{NumberedSubSection}

\begin{NumberedSubSection}{Smallness conditions on \aring/ homomorphism}
  Now we identify certain \ring/ homomorphisms which are particularly
  convenient to work with. See \ref{OrdinaryRingERegular} for the main
  motivating example.
  
\begin{defn}\label{DefineRegular}
  \Aring/ homomorphism $R\to k$ is \emph{\hregular/} if $k$ is \smalll/ as an
  $R$-module, \emph{\cohregular} if $k$ is \cosmall, and
  \emph{\ehregular/} if $k$ is \mixedsmall.
\end{defn}

\begin{rem}\label{DefineAbsoluteRegular}
  As in \ref{LatentMeansSameCells}, if $R\to k$ is either \hregular/
  or \cohregular/ it is also \ehregular. These are three very
  different conditions to put on the map $R\to k$, with \ehregularity/
  being by far the weakest one (see \ref{OrdinaryRingERegular}).

  Our notion of smallness is related to the notion of regularity 
  from commutative algebra.  For instance, a commutative
  ring $R$ is \emph{regular} (in the absolute sense) if and only if every
  finitely-generated discrete $R$-module $M$ is \smalll, i.e., has a
  finite length resolution by finitely generated projectives. Suppose
  that $f:R\to k$ is a surjection of commutative Noetherian rings. If
  $f$ is regular as a map of rings it is \hregular/ as a map of \rings,
  but the converse does not hold in general. The 
  point is that for $f$ to be regular in the ring-theoretic sense,
  certain additional conditions must be satisfied by the fibres of
  $R\to k$. 
\end{rem}

\begin{numbered}{Relationships between types of smallness}\label{TypesOfRegularity}
  Suppose that $k$ is an $R$-module and that $\EE=\EndR(k)$.  The
  \emph{double centralizer} of $R$ is the ring $\hat R=\EndE(k)$.
  Left multiplication gives a ring homomorphism $R\to\hat R$, and the
  pair $\Pair Rk$ is said to be \emph{\dccomplete/} if the
  homomorphism $R\to\hat R$ is an equivalence.  We show below
  (\ref{ThreeCompletionTypes}) that if $R\to k$ is a surjective map of
  Noetherian commutative rings with kernel $I\subset R$, then, as long
  as $k$ is a regular ring, $(R,k)$ is \dccomplete/ if and only if $R$
  is isomorphic to its $I$-adic completion.

  If $R$ is an augmented $k$-algebra, then $\EE=\EndR(k)$ is
  also an augmented $k$-algebra. The augmentation is provided by the
  natural map $\EndR(k)\to\EndOver k(k)\weq k$ induced by the
  $k$-algebra structure homomorphism $k\to R$.

\begin{prop}\label{RegularCriterion}
  Suppose that $R$ is \anaugmented/ $k$-algebra, and let
  $\EE=\EndR(k)$. Assume that the pair $(R,k)$ is \dccomplete. Then
  $R\to k$ is \hregular/ if and only if $\EE\to k$ is
  \cohregular. Similarly, $R\to k$ is \ehregular/ if and only if
  $\EE\to k$ is \ehregular. 
\end{prop}

\begin{titled}{Proof}
  If $k$ is finitely built from $R$ as an $R$-module, then by applying
  $\HomR(\whatever,k)$ to the construction process, we see that
  $\EE=\HomR(k,k)$ is finitely built from $k=\HomR(R,k)$ as an
  $\EE$-module. Conversely, if $\EE$ is finitely built from $k$ as an
  $\EE$-module, it follows that $k=\HomOver{\EE}(\EE, k)$ is
  finitely built from $\hat R\weq \HomOver{\EE}(k,k)$ as an
  $R$-module. If $R\weq\hat R$, this implies that $k$ is finitely
  built from $R$.

  For the rest, it is enough by symmetry to show that if $R\to k$ is
  \ehregular, then so is $\EE\to k$.
  Suppose then that $k$ is \mixedsmall/ over $R$ with Koszul complex
  $K$. Let $L=\HomR(K,k)$. Arguments as above show that $L$ is
  finitely built both from $\HomR(R,k)\weq k$ and from
  $\HomR(k,k)\weq\EE$ as an \EE-module. This means that $L$ will
  serve as a Koszul complex for $k$ over $\EE$, as long as $L$ builds
  $k$ over $\EE$. Let $\EEK=\EndR(K)$. By
  \ref{TLatentMoritaEquiv}, the natural map
  $L\TensorOver{\EEK}K\to k$ is an equivalence; it is
  evidently a map of \EE-modules. Since $\EEK$ builds $K$ over $\EEK$,
  $L\weq L\TensorOver{\EEK}\EEK$ builds $k$ over $\EE$.
  \qed
\end{titled}

In the following proposition, we think of $S\to R\to Q$ as a ``short
exact sequence'' of commutative \rings; often, such a sequence is obtained
by applying a cochain functor to a fibration sequence of spaces (cf.
\ref{EHRegularChains} or \ref{GorensteinCompactLie}).

  \begin{prop}\label{EhregularInTheMiddle}
    Suppose that $S\to R$ and $R\to k$ are homomorphisms of
    commutative \rings, and let $Q=R\TensorOver S k$. Note that $Q$
    is a commutative \ring/ and that there
    is a natural homomorphism $Q\to k$ which extends $R\to k$. 
    Assume that one of the following
    holds:
    \begin{enumerate}
    \item \label{MiddleOne} $S\to k$ is \ehregular/ and $Q\to k$ is \cohregular, or
    \item \label{MiddleTwo} $S\to k$ is \hregular/ and  $Q\to k$ is \ehregular.
    \end{enumerate}
    Then $R\to k$ is \ehregular.
  \end{prop}
  
  \begin{proof}
    Note that there is \aring/ homomorphism $R\to Q$.
    In case (\ref{MiddleOne}), suppose that $K$ is a Koszul complex for
    $k$ over $S$. We will show that $R\TensorOver S K$ is a Koszul
    complex for $k$ over $R$. Since $K$ is \smalll/ over $S$,
    $R\TensorOver SK$ is \smalll/ over $R$. Since $k$ finitely builds
    $K$ over $S$, $R\TensorOver Sk=Q$ finitely builds $R\TensorOver SK$ 
    over $R$.  But $k$ finitely builds $Q$ over $Q$, and hence
    over $R$; it follows that $k$ finitely builds $R\TensorOver SK$
    over $R$.  Finally, $K$ builds $k$ over $S$, and so $R\TensorOver S K$ 
   builds $Q$ over $R$; however, $Q$ clearly builds $k$ as a
    $Q$-module, and so {\it a fortiori} builds $k$ over $R$.

    In case (\ref{MiddleTwo}), let $K$ be a Koszul complex for $k$ over
    $Q$. We will show that $K$ is also a Koszul complex for $k$ over
    $R$. Note that $S\to k$ is \hregular, so that $k$ is \smalll/ over $S$
    and hence $Q=R\TensorOver Sk$ is \smalll/ over
    $R$. But $K$ is finitely built from
    $Q$ over $Q$ and hence over $R$; it follows that $K$ is \smalll/
    over $R$. Since $k$ finitely builds $K$ over $Q$, it does so over
    $R$; for a similar reason $K$ builds $k$ over $R$. 
  \end{proof}

\end{numbered}

\end{NumberedSubSection}

\begin{NumberedSubSection}{DC-completeness in algebra}
  Recall that a ring $k$ is said to be \emph{regular} if every
  finitely generated discrete module over $k$ has a finite projective
  resolution, i.e., if every finitely generated discrete module over $k$ is
  \smalll/ over $k$.

  \begin{prop}\label{ThreeCompletionTypes}
    Suppose that $R\to k$ is a surjection of commutative Noetherian
    rings with kernel ideal $I\subset R$. Assume that $k$ is a regular
    ring. Then the double centralizer map $R\to \hat R$ is can be
    identified with the $I$-adic completion map $R\to \lim_sR/I^s$.
    In particular, $\Pair Rk$ is \dccomplete/ (\ref{TypesOfRegularity}) if and only if
    $R$ is $I$-adically complete.
  \end{prop}

  \begin{proof}
    Let $\EE=\EndR(k)$, and $\hat R=\EndE(k)$, so that there is a
    natural homomorphism $R\to\hat R$ which is an equivalence if and
    only if $\Pair Rk$ is \dccomplete. We will show that $\hat R$ is
    equivalent to $\RIhat=\lim_sR/I^s$, so that $(R,k)$ is
    \dccomplete/ if and only if $R\to \RIhat$ is an isomorphism.

    Consider the class of all $R$-modules $X$ with the property that
    the natural map
    \begin{equation}\label{DoubleDual}
              X \to \HomE(\HomR(X,k),k)
    \end{equation}
    is an equivalence. The class includes $k$, and hence all
    $R$-modules finitely built from $k$. Each quotient $I^s/I^{s+1}$
    is finitely generated over $k$, hence \smalll/ over
    $k$, and hence finitely built from $k$ over $R$. It follows from
    an inductive argument that the modules $R/I^s$ are finitely built
    from $k$ over $R$, and consequently that \ref{DoubleDual} is an equivalence for
    $X=R/I^s$. By a theorem of Grothendieck \cite[2.8]{GrothendieckLocal}, there are isomorphisms
    \[
      \colim_s\Ext^i_R(R/I^s,k)\iso\begin{cases} k\quad i=0\\
                                                 0 \quad i>0
                                \end{cases}
    \]
    which (\ref{OurNotation})  assemble into an equivalence
    \[
       \hocolim_s\HomR(R/I^s,k)\weq\HomR(R,k)\weq k\,.
    \]
    This allows for the calculation
    \[
    \begin{aligned}
      \hat R\weq\HomE(k,k) &\weq \HomE(\hocolim_s\HomR(R/I^s,k),k)\\
                          & \weq \holim_s\HomE(\HomR(R/I^s,k),k)\\
                          & \weq \holim_s R/I^s\weq \RIhat\,\,.
    \end{aligned}
     \]
    It is easy to check that under this chain of equivalences the map
    $R\to \hat R$ corresponds to the completion map $R\to \RIhat$.
  \end{proof}
  
\end{NumberedSubSection}

\begin{NumberedSubSection}{DC-completeness in topology}\label{ChainCochainStuff}
   Suppose that $X$ is a connected pointed topological space, and that
  $k$ is a commutative \ring. For any space $Y$, pointed or not, let
  $Y_+$ denote $Y$ with a disjoint basepoint added, and  $\Stable Y_+$
  the associated suspension spectrum.
  We will consider two $k$-algebras associated to the pair
  $(X,k)$: the chain algebra $R=C_*(\Omega
  X;k)=k\TensorS \Stable(\Omega X)_+$ and the cochain algebra
  $S=C^*(X;k)=\Map_{\Sphere}(\Stable X_+,k)$. Here $\Omega X$ is the
  loop space on $X$, and $R$ is \aring/ because $\Omega
  X$ can be constructed as a topological or simplicial group;
  $R$ is essentially the group ring $k[\Omega X]$. The
  multiplication on $S$ is cup product coming from the diagonal
  map on $X$, and so $S$
  is a commutative $k$-algebra. Both of these objects are augmented,
  one by the map $R \to k$ induced by the map $\Omega
  X\to\operatorname{pt}$, the other by the map $S\to k$ induced
  by the basepoint inclusion $\operatorname{pt}\to X$. If $k$ is a
  ring, then $\Hof iR\iso H_i(\Omega X;k)$ and $\Hof iS\iso H^{-i}(X;k)$.

  The Rothenberg-Steenrod construction \cite{rothenberg-steenrod}
  shows that for any $X$ and $k$ there is an equivalence
  $S\weq\EndOver{R}(k)$. We will say that the
  pair $(X,k)$ is \emph{of Eilenberg-Moore type} if $k$ is a field,
  each homology group $H_i(X;k)$ is finite dimensional over $k$, and
  either
  \begin{enumerate}
  \item $X$ is simply connected, or
  \item $k$ is of characteristic $p$ and $\pi_1X$ is a finite $p$-group.
  \end{enumerate}
  If $(X,k)$ is of Eilenberg-Moore type, then by the Eilenberg-Moore
  spectral sequence construction (\cite{eilenberg-moore},
  \cite{dwyerSS}, \cite[Appendix C]{MandellPadic}), $R
  \weq\EndOver{S}(k)$ and both of the pairs $\Pair {R}k$
  and $\Pair{S} k$ are \dccomplete/ (\ref{TypesOfRegularity}).

  Keep in mind that if $(X,k)$ is of Eilenberg-Moore type, then the
  augmentation ideal of $\pi_0C_*(\Omega X;k)=k[\pi_1X]$ is nilpotent
  (cf. \ref{SpecialDownwardFiniteType}).
\end{NumberedSubSection}

\section{Examples of smallness}\label{CRegularityExamples}

In this section we look at some sample cases in which the smallness
conditions of \S\ref{CContext} are or are not satisfied.

\begin{NumberedSubSection}{Commutative rings}\label{OrdinaryRingERegular}
  If $R$ is a commutative Noetherian ring and $I\subset R$ is an
  ideal such that the quotient $R/I=k$ is a regular ring (\ref{DefineAbsoluteRegular}), then $R\to k$ is
  \ehregular/ \cite[\S6]{dwyer-greenlees}; the complex $K$ can be
  chosen to be the Koszul complex associated to any finite set of
  generators for $I$. The construction of the Koszul complex is
  sketched below in the proof of \ref{LocalCohomologySS}.  The pair $(R,k)$ is
  \dccomplete/ if and only if $R$ is complete and Hausdorff with
  respect to the $I$-adic topology (\ref{ThreeCompletionTypes}).

  For example, if $R$ is a Noetherian local ring with residue field
  $k$, then the map $R\to k$ is \ehregular; this map is \hregular/ if
  and only if $R$ is a regular ring (Serre's Theorem) and \cohregular/ if and
  only if $R$ is artinian.
\end{NumberedSubSection}

\begin{NumberedSubSection}{The sphere spectrum}
  Consider the map $\Sphere\to\Fp$ of commutative \rings; here as
  usual \Sphere/ is the sphere spectrum and the ring \Fp/ is
  identified with the associated Eilenberg-MacLane spectrum. This map
  is \emph{not} \ehregular. A Koszul complex $K$ for $\Sphere\to\Fp$
  would be a stable finite complex with nontrivial mod~$p$ homology
  (because $K$ would build $\Fp$), and only a finite number of
  non-trivial homotopy groups, each one a finite $p$-group
  (because $\Fp$ would finitely build $K$). We leave it to the reader
  to show that no such $K$ exists, for instance because of Lin's theorem \cite{lin}
  that $\Map_{\Sphere}(\Fp,\Sphere)\weq 0$. 

  Let $\Sphere_p$ denote the $p$-completion of the sphere
  spectrum. The map $\Sphere\to\Fp$ is not \dccomplete/, but
  $\Sphere_p\to\Fp$ is; this can be interepreted in terms of the
  convergence of the classical mod~$p$ Adams spectral sequence. 
\end{NumberedSubSection}

The next two examples refer to the following proposition.

  \begin{prop}\label{FiniteSkeletaFiniteType}\label{FiniteSkeletaFiniteStableType}
  \label{FiniteComplexCoregular}
  Suppose that $X$ is a pointed connected finite complex, that $k$ is
  a commutative \ring, that $R$ is the augmented $k$-algebra
  $C_*(\Omega X;k)$, and that $S$ is the augmented $k$-algebra
  $C^*(X;k)=\End_R(k)$.  Then $k$ is small as an $R$-module and cosmall
  as an $S$-module.
\end{prop}

\begin{proof}
  Let $E$ be the total space of the universal principal bundle over
  $X$ with fibre $\Omega X$, so that $E$ is contractible and
  $M=C_*(E;k)\weq k$. The action of $\Omega X$ on $E$ induces an
  action of $R$ on $M$ which amounts to the augmentation action of $R$
  on $k$. Let $E_i$ be the inverse image in $E$ of the $i$-skeleton of
  $X$, and let $M_i$ be the $R$-module $C_*(E_i;k)$. Then
  $M_i/M_{i-1}$ is equivalent to a finite sum $\oplus \Shift^iR$
  indexed by the $i$-cells of $X$. Since $k\weq M= M_n$ (where $n$ is
  the dimension of~$X$), it
  follows that $k$ is small as an $R$-module. The last statement is
  immediate, as in the proof of \ref{SpecialDownwardFiniteType}.
\end{proof}

\begin{rem}\label{FiniteTypeFiniteUpward}
  The argument above also shows that if $X$ is merely of finite type
  (i.e., has a finite number of cells in each dimension), then the
  augmentation module $k$ is
  of upward finite type over $R=C_*(\Omega X;k)$. 
\end{rem}

\begin{NumberedSubSection}{Cochains}\label{RegularCochainAlgebras}
  Suppose that $X$ is a pointed connected topological space, that $k$
  is a commutative \ring, and that
  $R$ is the augmented $k$-algebra $C^*(X;k)$.  
  \begin{enumerate}
  \item The map $R\to k$ is \cohregular/ if $X$ is a
    finite complex (\ref{FiniteComplexCoregular}).
  \item If $k$ is a field, then $R\to k$ is \cohregular/ if and only
    if $H^*(X;k)$ is finite-dimensional (\ref{BestOfAllCoregularity}).
  \item If $(X,k)$ is of Eilenberg-Moore type, then $R\to k$ is
    \hregular/ if and only if $H_*(\Omega X;k)$ is finite-dimensional
    (\ref{BestOfAllCoregularity}, \ref{RegularCriterion}).
  \end{enumerate}
\end{NumberedSubSection}

\begin{NumberedSubSection}{Chains}\label{RegularChainAlgebras}
  Suppose that $X$ is a pointed connected topological space, that $k$
  is a commutative \ring, and that
  $R$ is the augmented $k$-algebra $C_*(\Omega X;k)$.  
  \begin{enumerate}
  \item The map $R\to k$ is \hregular/ if $X$ is a finite complex
    (\ref{FiniteSkeletaFiniteType}).
  \item If $(X,k)$ is of Eilenberg-Moore type, then $R\to k$ is
    \hregular/ if and only if $H^*(X;k)$ is finite-dimensional
    (\ref{RegularCriterion}, \ref{RegularCochainAlgebras}).
  \item If $(X,k)$ is of Eilenberg-Moore type, then $R\to k$ is
    \cohregular/ if and only if $H_*(\Omega X;k)$ is
    finite-dimensional (\ref{BestOfAllCoregularity}).
  \end{enumerate}
  The parallels between
  \ref{RegularCochainAlgebras} and \ref{RegularChainAlgebras} are
  explained by \ref{RegularCriterion}.

\end{NumberedSubSection}

\begin{NumberedSubSection}{Completed classifying spaces}\label{EHRegularChains}
  Suppose that $G$ is a compact Lie group (e.g., a finite group), that
  $k=\Fp$, and that $X$ is the $p$-completion of the classifying space
  $BG$
  in the sense of Bousfield-Kan \cite{BK}. Let 
  \[ R=C^*(X;k) \text{ and } \EE=C_*(\Omega X;k)\,.\] We will show in
  the following paragraph that $R\to k$ and $\EE\to k$ are both
  \ehregular, and that the pair $(X,k)$ is of Eilenberg-Moore type.
  There are many $G$ for which neither $H_*(\Omega X;k)$ nor
  $H^*(X;k)$ is finite dimensional \cite{levi}; by
  \ref{RegularCochainAlgebras} and \ref{RegularChainAlgebras}, in such
  cases the maps  $R\to k$ and $\EE\to k$ are neither 
  \hregular/ nor \cohregular. We are interested in these examples for
  the sake of local cohomology theorems (\ref{GorensteinCompactLie}).

  By elementary representation theory there is a faithful embedding
  $\rho:G\to SU(n)$ for some $n$, where $SU(n)$ is the special unitary
  group of $n\times n$ Hermitian matrices of determinant one. Consider
  the associated fibration sequence
  \begin{equation}\label{EmbeddingFibration}
      M= SU(n)/ G\to BG\to BSU(n)\,. 
  \end{equation}
  The fibre $M$ is a finite complex. Recall that $R=C^*(BG;k)$; write
  $S=C^*(BSU(n);k)$ and $Q=C^*(M;k)$.  Since $BSU(n)$ is
  simply-connected, the Eilenberg-Moore spectral sequence of
  \ref{EmbeddingFibration} converges and $Q\weq k\TensorOver S R$ (cf.
  \cite[5.2]{MandellPadic}).  The map $S\to k$ is \hregular/ by
  \ref{RegularCochainAlgebras}(3) and $Q\to k$ is \cohregular/ by
  \ref{RegularCochainAlgebras}(2);  
  it follows from \ref{EhregularInTheMiddle}
  that $R\to k$ is \ehregular.   Since $\pi_1BG=\pi_0G$ is
  finite, $BG$ is $\Fp$-good (i.e., $ C^*(X;k)\weq C^*(BG;k)$), and $\pi_1X$
  is a finite $p$-group \cite[VII.5]{BK}. In particular, $(X,k)$ is of
  Eilenberg-Moore type. Since $\EE=C_*(\Omega
  X;k)$ is thus equivalent to $\EndR(k)$,  we conclude from \ref{RegularCriterion} that $\EE\to k$ is
  also \ehregular.
\end{NumberedSubSection}

\begin{NumberedSubSection}{Group rings}
  If $G$ is a finite group and $k$ is a commutative ring, then the
  augmentation map $k[G]\to k$ is \ehregular. We will prove this by
  producing a Koszul complex $K$ for $\Z$ over $\Z[G]$; it is then
  easy to argue that $k\TensorOver{\Z}K$ is a Koszul complex for $k$
  over $k[G]$. Embed $G$ as above into a unitary group $SU(n)$ and let
  $K=C_*(SU(n);\Z)$.  The space $SU(n)$ with the induced left
  $G$-action is a compact manifold on which $G$ acts smoothly and
  freely, and so by transformation group theory \cite{illman} can be constructed
  from a finite number of $G$-cells of the form $(G\times D^i,G\times
  S^{i-1})$. This implies that $K$ is small over $\Z[G]$, since, up
  to equivalence over $\Z[G]$, $K$ can be identified with the $G$-cellular
  chains on $SU(n)$. Note that $G$ acts trivially on
  $\Hof *K=H_*(SU(n);\Z)$ (because $SU(n)$ is connected) and that,
  since $H_*(SU(n);\Z)$ is torsion free, each group $\Hof iK$ is
  isomorphic over $G$ to a finite direct sum of copies of the
  augmentation module~$\Z$. The Postnikov argument in the proof of
  \ref{BestOfAllCoregularity} thus shows that $K$ is finitely built
  from $\Z$ over $\Z[G]$. Finally, $K$ itself is \aring, the action of
  $\Z[G]$ on $K$ is induced by a homomorphism $\Z[G]\to K$,
  and the augmentation $\Z[G]\to \Z$
  extends to an augmentation $K\to \Z$. Since $K$ builds $\Z$ over
  $K$ (see \ref{FinitelyBuilt}), $K$ certainly builds $\Z$ over $\Z[G]$.
\end{NumberedSubSection}

\section{Matlis lifts}\label{CMatlisDuality}

Suppose that $R$ is a commutative Noetherian local ring, and that
$R\to k$ is reduction modulo the maximal ideal. Let $\II(k)$ be the
injective hull of $k$ as an $R$-module. The starting point of this
section is the isomorphism 
\begin{equation}\label{AlgebraicLiftFormula}
\Homk(X,k)\weq\HomR(X,\II(k))\,,
\end{equation}
which holds for any $k$-module $X$. We think of $\II(k)$ as a lift of
$k$ to an $R$-module, not the obvious lift obtained by using the
homomorphism $R\to k$, but a more mysterious construction that allows
for \ref{AlgebraicLiftFormula}.  If we define the \emph{Pontriagin dual} of an
$R$-module $M$ to be $\HomR(M,\II(k))$,
then Pontriagin duality is a correspondingly mysterious construction for
$R$-modules which extends ordinary $k$-duality for $k$-modules.

We generalize this in the following way.
\begin{defn}\label{ReallyDefineMatlisLift}
  Suppose that $R\to k$ is a map of \rings, and that $N$ is a $k$-module.
  An $R$-module \II/ is said to be a \emph{Matlis lift} of $N$ if the
  following two conditions hold:
\begin{enumerate}
\item $\HomR(k,\II)$ is equivalent to $N$ as a left $k$-module, and
\item $\II$ is effectively constructible from~$k$.
\end{enumerate}
\end{defn}

\begin{rem}\label{ExtendedDuality}
  If \II/ is a Matlis lift of~$N$ and $X$ is an arbitrary $k$-module,
  then the adjunction equivalence
  \[\HomR(X,\II)\weq\Homk(X,\HomR(k,\II))\]
  implies, by \ref{ReallyDefineMatlisLift}(1), that there is an
  equivalence
\begin{equation}\label{LiftFormula}
  \Homk(X,N) \weq \HomR(X,\II)\,.
\end{equation}
This is the crucial property of a Matlis lift (cf.
\ref{AlgebraicLiftFormula}).  Condition
\ref{ReallyDefineMatlisLift}(2) tightens things up a bit. There's no
real reason not to assume that \II/ is $k$-cellular, since if \II/ satisfies
\ref{ReallyDefineMatlisLift}(1), so does $\Cell_k\II$. The somewhat stronger assumption that \II/ is effectively
constructible from~$k$ 
will allow Matlis lifts to be constructed and enumerated
(\ref{ClassifyMatlisLifts}).  In many situations, the assumption that
\II/ is $k$-cellular implies \ref{ReallyDefineMatlisLift}(2).
\end{rem}

\begin{rem}
  If $R$ is a commutative Noetherian local ring, and $R\to k$ is
  reduction modulo the maximal ideal, then the injective hull $\II(k)$
  is a Matlis lift of $k$ (\ref{LocalRingMatlisLift}).
\end{rem}

For the rest of the section we assume that $R\to k$ is a map of
\rings, and that $N$ is a $k$-module. Let $\EE=\EndR(k)$. Observe that
the right multiplication action of $k$ on itself gives a homomorphism
$k\op\to\EE$, or equivalently $k\to \EE\op$, so it makes sense to look
at right \EE-actions on $N$ which extend the left $k$-action.

\begin{defn}\label{DefineMatlisLift}
   An \emph{\Elift} of $N$ is a right
  $\EE$-module structure on $N$ which extends the left $k$-action. An
  \Elift/ of $N$ is said to be \emph{of Matlis type} if
  the natural map
  \begin{equation}\label{MatlisLiftFormula}
   N\weq N\TensorE\HomR(k,k)\to\HomR(k, N\TensorE k)
  \end{equation}
  is an equivalence. (Here the action of $R$ on
  $N\TensorE k$ is obtained from the left action of $R$ on $k$.) 
\end{defn}

\begin{rem}\label{DefineMatlisType}
  More generally,  an arbitrary right $\EE$-module $N$ is said to be of \emph{Matlis type} if the map
  \ref{MatlisLiftFormula} is an equivalence.
\end{rem}

The following proposition gives a classification of Matlis lifts.

\begin{prop}\label{ClassifyMatlisLifts}
  The correspondences
  \[
      \II\mapsto\HomR(k,\II)  \quad\quad\quad N \mapsto N\TensorE k
  \]
  give inverse bijections, up to equivalence, between Matlis lifts
  \II/ of $N$ and \Elifts/ of $N$ which are of Matlis type.
\end{prop}

\begin{proof}
  If $\II$ is a Matlis lift of $N$, then $\HomR(k,\II)$ is equivalent
  to $N$ as a $k$-module (\ref{ReallyDefineMatlisLift}), and so the
  natural right action of \EE/ on $\HomR(k,\II)$ provides an \Elift/
  of $N$. Equivalence \ref{MatlisLiftFormula} holds because \II/
  is effectively constructible from~$k$; consequently, this \Elift/ is
  of Matlis type.

  Conversely, given an \Elift/ of $N$ which is of Matlis type, let
  $\II=N\TensorE k$. Equivalence \ref{MatlisLiftFormula} guarantees
  that $\II$ satisfies \ref{ReallyDefineMatlisLift}(1), and the same
  formula leads to the conclusion that $\II$ is effectively
  constructible from~$k$.
\end{proof}

The following observation is useful for recognizing Matlis lifts.

\begin{prop}\label{MatlisTypeCellularApprox}
  Suppose that $R\to k$ is a map of \rings, that $\EE=\EndR(k)$, and
  that $M$ is an $R$-module. Then the right $\EE$-module
  $\HomR(k,M)$ is of Matlis type if and only if the evaluation map
  \[
     \HomR(k,M)\TensorE k\to M
  \]
  is a $k$-cellular approximation, i.e., if and only if $\CellOf kM$
  is effectively constructible from~$k$.   
\end{prop}

\begin{proof}
  Let $N=\HomR(k,M)$. Since $N$ is $\EE\op$-cellular over $\EE\op$,
  $N\TensorE k$ is $k$-cellular over $R$. This implies that
  the evaluation map $\epsilon$ is a $k$-cellular approximation if and only if it
  is a $k$-equivalence. Consider the maps
  \[
    N\TensorE\HomR(k,k)\to\HomR(k,N\TensorE
    k)\RightArrow{\HomR(k,\epsilon)} N\,.
  \]
  It is easy to check that the composite is the obvious equivalence,
  so the left hand map is an equivalence ($N$ is of Matlis type) if
  and only if the right-hand map is an equivalence ($\epsilon$ is a
  $k$-equivalence).  The final statement is from \ref{ThreeConstructible}.
\end{proof}

\begin{NumberedSubSection}{Matlis duality}\label{DefineMatlisDuality}
  In the situation of \ref{DefineMatlisLift}, let $N=k$ and let
  $\II=k\TensorE k$ be a Matlis lift of $k$. The \emph{Pontriagin
    dual} or \emph{Matlis dual} of an $R$-module $M$ (with respect to
  $\II$) is defined to be $\HomR(M,\II)$. By \ref{ExtendedDuality},
  Matlis duality is a construction for $R$-modules which extends
  ordinary $k$-duality for $k$-modules. Note, however, that in the
  absence of additional structure (e.g., commutativity of $R$) it is
  not clear that $\HomR(M,\II)$ is a right $R$-module. We will come up with
  one way to remedy this later on (\ref{GorensteinMeansNiceDuals}).
\end{NumberedSubSection}

\begin{NumberedSubSection}{Existence of Matlis lifts}\label{MatlisLiftExistence}
  We give four conditions under which a right $\EE$-module is of
  Matlis type (\ref{DefineMatlisType}), and so gives rise to a Matlis lift of the underlying
  $k$-module.  The first two conditions are of an algebraic nature;
  the second two may seem technical, but they apply to many ring
  spectra, chain algebras, and cochain algebras. In all of the
  statements below, $R\to k$ is a map of \rings, $\EE=\EndR(k)$, and
  $N$ is a right \EE-module.

  \begin{prop}\label{PRegularLift}
    If $R\to k$ is \hregular, then any $N$ is of Matlis type.
  \end{prop}

  \begin{proof}
    Calculate
   \[ \HomR(k,N\TensorE k)\weq N\TensorE\HomR(k,k)\weq N\TensorE
   \EE\weq N\,
   \] 
   where the first weak equivalence comes from the fact that $k$ is
   small as an $R$-module.
  \end{proof}

  \begin{prop}\label{EhregularMatlisCondition}
    If $R\to k$ is \ehregular/, then $N$ is of Matlis type
    if and only if there exists an $R$-module $M$ such that 
    $N$ is equivalent to $\HomR(k,M)$ as a
    right \EE-module. 
  \end{prop}

  \begin{proof}
    If $N$ is of Matlis type, then $M=N\TensorE k$ will do.
    Given $M$, the fact that $\HomR(k,M)$ is of Matlis type follows from
    \ref{MatlisTypeCellularApprox} and \ref{TLatentMoritaEquiv} (cf. \ref{ThreeConstructible}).
  \end{proof}

  \begin{prop}\label{PChainCanonical}
    Suppose that $k$ and $N$
     are bounded above, that $k$ is of upward finite type as an
    $R$-module, and that $N$ is of downward type as an
    $\EE\op$-module. Then $N$ is of Matlis type.
  \end{prop}

  \begin{prop}\label{PCochainCanonical}
    Suppose that $k$ and $N$ are
    bounded below, that $k$ is of downward finite type as an
    $R$-module and that $N$ is of upward type as an
    $\EE\op$--module.  Then $N$ is of Matlis type.
  \end{prop}

    \begin{titled}{Proof of \ref{PChainCanonical}}
      Note first that since $k$ and $N$ are bounded above and $N$ is
      of downward type as an $\EE\op$-module, $N\TensorE k$ is also
      bounded above.  Consider the class of all $R$-modules $X$ such
      that the natural map
    \[
      N\TensorE \HomR (X, k)\to \HomR(X, N\TensorE k)
    \]
    is an equivalence. This certainly includes $R$, and so by
    triangle arguments includes everything that can be finitely built from
    $R$. We must show
    that the class contains $k$. Pick an integer $B$, and suppose that
    $A$ is another integer. Since $k$ is of upward finite type as
    an $R$-module and both $k$ and $N\TensorE k$ are bounded above,
    there exists an $R$-module $X$, finitely built from $R$, and a map
    $X\to k$ which induces isomorphisms
    \begin{equation}\label{VariousIsos}
     \begin{aligned}
        \Hof i\HomR(k, k) &\RightArrow\iso \Hof i\HomR(X,k)\\
        \Hof i\HomR(k,N\TensorE k) &\RightArrow\iso \Hof
        i\HomR(X,N\TensorE k)
     \end{aligned}\quad \quad i>A.
    \end{equation}
    Now $N$ is of downward type as a right
    \EE-module, so if we choose $A$ small enough we can guarantee that
    the map
    \[
       \Hof i (N\TensorE \HomR(k,k))  \to  \Hof i (N\TensorE \HomR(X,k))
    \]
    is an isomorphism for $i>B$. By reducing $A$ if necessary
    (which of course affects the choice of $X$), we can assume $A\le
    B$.  Now consider the commutative diagram
    \begin{equation}\label{StableCommDiagram}
    \begin{CD}
       N\TensorE \HomR(k,k) @>>> \HomR(k, N\TensorE k)\\
          @VVV                      @VVV              \\
       N\TensorE \HomR(X,k) @>>> \HomR(X, N\TensorE k)
    \end{CD}
    \end{equation}
    The lower arrow is an equivalence, because $X$ is finitely built
    from $R$, and the vertical arrows are isomorphisms on $\Hof i$ for
    $i>B$. Since $B$ is arbitrary, it follows that the upper arrow is
    an equivalence.\qed
  \end{titled}
  
  \begin{titled}{Proof of \ref{PCochainCanonical}}
    This is very similar to the proof above, but with the inequalities
    reversed. Observe that since $k$ and $N$ are bounded below, and
    $N$ is of upward type as an $\EE$-module, $N\TensorE k$ is also
    bounded below. Pick an integer $B$, and let $A$ be another
    integer. Since $k$ is of downward finite type as an $R$-module and
    both $k$ and $N\TensorE k$ are bounded below, there exists an $X$
    finitely built from $R$ such that the maps in \ref{VariousIsos} are
    isomorphisms for $i<A$. Now $N$ is of upward type as a right
    \EE-module, so if we choose $A$ large enough we can guarantee that
    the map
    \[
       \Hof i (N\TensorE \HomR(k,k))  \to  \Hof i (N\TensorE \HomR(X,k))
    \]
    is an isomorphism for $i<B$. By making $A$ larger if necessary, we
    can assume $A>B$. The proof is now completed by using the
    commutative diagram \ref{StableCommDiagram}. \qed
  \end{titled}

\end{NumberedSubSection}
\section{Examples of Matlis lifting}\label{MatlisLiftingExamples}

In this section we look at particular examples of Matlis lifting
(\S\ref{CMatlisDuality}). In each case we start with a morphism $R\to
k$ of rings, and look for Matlis lifts of $k$.  As usual, $\EE$
denotes $\EndR(k)$.

\begin{NumberedSubSection}{Local rings}\label{LocalRingMatlisLift}
  Suppose that $R$ is a commutative Noetherian local ring with maximal
  ideal $I$ and residue field $R/I=k$, and that $R\to k$ is the
  quotient map. Let $\II=\II(k)$ be the injective hull of $k$ (as an
  $R$-module). We will show that $\II$ is  a Matlis lift of $k$.

  To see this, first note that $\II$ is $k$-cellular, or
  equivalently \cite[6.12]{dwyer-greenlees}, that each element of
  $\II$ is annihilated by some power of $I$. Pick an element
  $x\in\II$; by Krull's Theorem \cite[10.20]{atiyah-macdonald} the
  intersection $\cap_jI^jx$ is trivial. But each submodule $I^jx$ of
  $\II$ is either trivial itself or contains $k\subset \II$
  \cite[p.~281]{matsumura}. The conclusion is that $I^jx=0$ for
  $j>\!\!>0$. Since $\HomR(k,\II)\weq k$ (again, for instance, by
  \cite{matsumura}) and $\II$ is effectively constructible from~$k$ 
  (\ref{OrdinaryRingERegular}, \ref{TLatentMorita}), $\II$ provides an \Elift/ of $k$.
  Up to equivalence there is exactly one \Elift/ of $k$
  (\ref{UniqueModuleStructures}), and so in fact $\II$ is the only
  Matlis lift of $k$.

  For instance, if $R\to k$ is $\Z_{(p)}\to\Fp$, then $\II\weq
  k\TensorE k$ is $\Zpinfty$ (cf. \cite[\S3]{dwyer-greenlees}), and
  Matlis duality (\ref{DefineMatlisDuality}) for $R$-modules is
  Pontriagin duality for $p$-local abelian groups.
\end{NumberedSubSection}

\begin{NumberedSubSection}{$k$-algebras}\label{kAlgebraExample}
  Suppose that $R$ is \anaugmented/ $k$-algebra, and let $M$ be the
  $R$-module $\Homk(R,k)$. The left $R$-action on $M$ is induced by
  the right $R$-action of $R$ on itself.  By an adjointness calculation,
  $\HomR(k,M)$ is equivalent to $k$, and so in this way $M$ provides
  an \Elift/ of $k$. If this \Elift/ is of Matlis type, then the
  $R$-module $k\TensorE k$, which by \ref{MatlisTypeCellularApprox} is
  equivalent to $\CellOf k^R\Homk(R,k)$, is a Matlis lift of $k$. 
  There are equivalences
  \[
    \Homk(k\TensorE k, k)\weq \HomE(k, \Homk(k,k))\weq\HomE(k,k)\weq
    \hat R\,,
  \]
  so that if $(R,k)$ is \dccomplete, the Matlis lift $k\TensorE k$
  is pre-dual over $k$ to $R$ (i.e., the $k$-dual of $k\TensorE k$
  is~$R$). Note that
  this calculation does not depend on assuming that $R$ is small in
  any sense as a $k$-module; there is an interesting example below in
  \ref{MatlisSuspensionLoops}.
\end{NumberedSubSection}

\begin{NumberedSubSection}{The sphere spectrum}\label{BrownComenetz}
  Let $R\to k$ be the unit map $\Sphere\to\Fp$. (Recall that we
  identify the ring $\Fp$ with the corresponding Eilenberg-MacLane
  \ring.)  The endomorphism \ring/ \EE/ is the Steenrod
  algebra spectrum, with $\Hof{-i}E$ isomorphic to the degree~$i$
  homogeneous component of the Steenrod algebra. Since $k$ has a
  unique \Elift/ (\ref{UniqueModuleStructures}) and the conditions of
  \ref{PChainCanonical} are satisfied (\ref{UpwardFiniteType},
  \ref{DownwardFiniteType}), $k$ has a unique Matlis lift given by
  $k\TensorE k$. Let $J$ be the Brown-Comenetz dual of $\Sphere$
  \cite{brown-comenetz;duality} and $J_p$ its $p$-primary summand.  We
  argue below that $J_p$ is $k$-cellular; by the basic property of
  Brown-Comenetz duality, $\HomR(k,J_p)\weq k$.  By
  \ref{MatlisTypeCellularApprox} the evaluation map $k\TensorE k\to
  J_p$ is a $k$-cellular approximation and hence, because $J_p$ is
  $k$-cellular, an equivalence.  Matlis duality amounts to the
  $p$-primary part of Brown-Comenetz duality.  Arguments parallel to 
  those in the proof of \ref{PChainCanonical} show  
  that if $X$ is spectrum which is bounded below and of finite type 
  then the natural map
  \[
   k\TensorE\HomR(X,k)\to  \HomR(X, k\TensorE k)
  \] 
  is an equivalence. Suppose that $X_*$ is an Adams resolution of the
  sphere.   Taking the 
  Brown-Comenetz dual $\HomR(X_*,k\TensorE k)$
  gives a spectral sequence which is the $\Fp$-dual of the
  mod~$p$ Adams spectral sequence. On the other hand, computing
  $\Hof*\HomR(X_*, k)$ amounts to taking the cohomology
  of $X_*$ and so gives a free resolution of $k$ over
  the Steenrod Algebra; the spectral sequence associated to
  $k\TensorE\HomR(X_*,k)$ is
  then
  the Kunneth 
  spectral sequence
  \[
    \Tor_*^{\Hof*\EE}(\Hof*k,\Hof*k) \Rightarrow \Hof*(k\TensorE k)\iso\Hof*J_p\,.
  \]
  It follows that these two spectral sequences are isomorphic.

  To see that $J_p$ is $k$-cellular, write $J_p=\hocolim J_p(-i)$, where
  $J_p(-i)$ is the $(-i)$-connective cover of $J_p$. Each $J_p(-i)$ has
  only a finite number of homotopy groups, each of which is a
  finite $p$-primary torsion group, and it follows immediately that
  $J_p$ can be finitely built from $k$. Thus $J_p$, as a homotopy
  colimit of $k$-cellular objects, is itself $k$-cellular (cf. \ref{BuiltFromk}).   
\end{NumberedSubSection}

\begin{NumberedSubSection}{Cochains}\label{CochainAlgebraExample}
  Suppose that $X$ is a pointed connected space and $k$ is a field
  such that $(X,k)$ is of Eilenberg-Moore type (\ref{ChainCochainStuff}).
  Let $R=C^*(X;k)$ and $\EE=\EndR(k)\weq C_*(\Omega X;k)$, and suppose that some \Elift/
  of $k$ is given. By \ref{UpwardFiniteType}, $k$ is of upward finite type
  over $\EE\op$, and hence 
  of downward finite type over
  $R$ (\ref{SpecialDownwardFiniteType}), the conditions of
  \ref{PCochainCanonical} are satisfied, and $\II=k\TensorE k$ is a
  Matlis lift of $k$. In fact there is only one \Elift/ of~$k$; this
  follows from
  \ref{UniqueModuleStructures}
  and the fact 
  that if $k$ is  a field of characteristic $p$ and $G$ is a finite
  $p$-group, any homomorphism $G\to k^\times$ is trivial.
  In these cases
  the Matlis lift $\II=k\TensorE k$ is equivalent by the
  Rothenberg-Steenrod construction to $C_*(X;k)=\Homk(R,k)$. Observe
  in particular that $\Homk(R,k)$ is $k$-cellular as an $R$-module;
  this also follows from \ref{BuiltFromk}.

\end{NumberedSubSection}

\begin{NumberedSubSection}{Chains}\label{ChainAlgebraExamples}
  Let $X$ be a connected pointed space, $k$ a field, and $R$ the chain algebra
  $C_*(\Omega X;k)$, so that $\EE\weq C^*(X;k)$. By
  \ref{UniqueModuleStructures} there is only one \Elift/ of $k$,
  necessarily given by the augmentation action of $\EE$ on $k$.
  Suppose that $k$ has upward finite type as an $R$-module, for
  instance, suppose that the conditions of \ref{UpwardFiniteType}
  hold, or that $X$ has finite skeleta
  (\ref{FiniteTypeFiniteUpward}). Then, by \ref{DownwardFiniteType}
  and \ref{PChainCanonical}, $k$ has a unique Matlis lift, given by
  $k\TensorE k$, or alternatively (\ref{kAlgebraExample}) by $\CellOf
  k\Homk(R,k)\weq\CellOf kC^*(\Omega X;k)$. We have \emph{not} assumed
  that $(X,k)$ is of Eilenberg-Moore type, and so the identification
  \[
      k\TensorE k \weq \CellOf kC^*(\Omega X;k)
  \]
  gives an interpretation of the abutment of the cohomology
  Eilenberg-Moore spectral sequence associated to the path fibration
  over $X$; this is in some sense dual to the interpretation of the
  abutment of the corresponding homology spectral as a suitable
  completion of $C_*(\Omega X)$ \cite{exotic-convergence}.
\end{NumberedSubSection}

\begin{NumberedSubSection}{Suspension spectra of loop spaces}\label{MatlisSuspensionLoops}
  Suppose that $X$ is a pointed finite complex, let $k=\Sphere$, and
  let $R$ be the augmented $k$-algebra $C_*(\Omega X;k)$. Then $\EE$
  is equivalent to $C^*(X;k)$, i.e., to the Spanier-Whitehead dual of
  $X$ (\ref{ChainCochainStuff}). Since $X$ is finite, $k$ is small as
  an $R$-module (\ref{FiniteSkeletaFiniteType}).  It follows from
  \ref{PRegularLift} that Matlis lifts of $k$ correspond bijectively
  to \Elifts/ of $k$.  Note that since the augmentation action of
  $\EE$ on $k$ factors through $\EE\to k$, and $k$ is commutative,
  this augmentation action amounts in itself to an \Elift. (It is
  possible to show that this is the only \Elift/ of $k$, but we will
  not do that here.) By inspection, this augmentation \Elift/ of $k$
  is the same as the \Elift/ obtained by letting \EE/ act in the
  natural way on $\HomR(k,\Hom_{k}(R,k))\weq k$ as in
  \ref{kAlgebraExample}. By \ref{MatlisTypeCellularApprox}, the
  corresponding Matlis lift $k\TensorE k$ is $\CellOf{k}\HomOver{k}(R,k)$.

  Suppose in addition that $X$ is $1$-connected, and write $k\TensorE
  k$ as the realization of the ordinary simplicial bar construction
  \[
     k\TensorS k \Leftarrow k\TensorS \EE\TensorS k\Lleftarrow
     k\TensorS\EE\TensorS\EE\TensorS k 
    \cdots \,.
  \]
  The spectrum $\HomS(k\TensorS k, \Sphere)$ is then the total complex
  of the corresponding cosimplicial object
  \[
     \HomS(k\TensorS k, \Sphere)\Rightarrow \HomS(k\TensorS\EE\TensorS
     k, \Sphere)\Rrightarrow\cdots\,.
  \]
  This is the cosimplicial object obtained by applying the unpointed
  suspension spectrum functor to the cobar construction on $X$, and by
  a theorem of Bousfield \cite{bousfieldss} its total complex is the
  suspension spectrum of $\Omega X$, i.e., $R$. Equivalently,
  Bousfield's theorem shows that in this case $(R,k)$ is \dccomplete.
  In this way if $X$ is $1$-connected the  Matlis lift of $k$ is 
  a Spanier-Whitehead pre-dual of $R$ (cf. \ref{kAlgebraExample}).
  This object has come up in a
  different way in work of N.  Kuhn \cite{kuhn}.

\end{NumberedSubSection}
\section{Gorenstein $\mathbb{S}$-algebras}\label{CGorenstein}

If $R$ is a commutative Noetherian local ring with maximal ideal $I$
and residue field $R/I=k$, one says that $R$ is \emph{Gorenstein} if
$\Ext_R^*(k,R)$ is concentrated in a single degree, and is isomorphic
to $k$ there. We give a similar definition for \rings, with an extra
technical condition added on.

\begin{defn}\label{DefineGorenstein}
  Suppose that $R\to k$ is a map of \rings, and let $\EE=\EndR(k)$.
  Then $R\to k$ is \emph{\hGorenstein/} of shift $\SH$ if the
  following two conditions hold:
  \begin{enumerate}
  \item as a left $k$-module, $\HomR(k,R)$ is equivalent to
    $\Shift^\SH k$, and
  \item \label{CanonicalRequirement} as a right $\EE$-module, 
    $\HomR(k,R)$ is of Matlis type (\ref{DefineMatlisType}).
  \end{enumerate}
\end{defn}

\begin{rem}\label{GorensteinMeansCells}
  Suppose that $R\to k$ is \hGorenstein/ of shift $\SH$, and give
  $\Shift^ak$ the right $\EE$-module structure from
  \ref{DefineGorenstein}(1). Then by \ref{MatlisTypeCellularApprox},
  $\CellOf k(R)$ is equivalent to $\Shift^ak\TensorE k$. It follows
  from  \ref{SpecialGorensteinConditions} and
  \ref{OrdinaryRingERegular} that if $R\to k$ is a map from a
  commutative Noetherian
  local ring to its residue field, then $R\to k$ is \hGorenstein/ in
  our sense if and only if $R$ is Gorenstein in the classical sense of
  commutative algebra.
\end{rem}

\begin{rem}\label{GorensteinMeansNiceDuals}
  Definition \ref{DefineGorenstein} does not exhaust all of the
  structure in $\HomR(k,R)$; in fact, the right action of $R$ on
  itself gives a right $R$-action on $\HomR(k,R)$ which commutes with
  the right $\EE$-action (since $\EE$ acts through $k$). This implies
  that if $R\to k$ is \hGorenstein/ and $k$ is given the right
  $\EE$-action obtained from $k\weq \Shift^{-a}\HomR(k,R)$, then the
  Matlis lift $\II= k\TensorE k$ of $k$ inherits a right
  $R$-action. In this case the Matlis dual
  $\HomR(M,\II)$ of a left $R$-module is naturally a right $R$-module. 
\end{rem}

In the \ehregular/ case it is possible to simplify definition
\ref{DefineGorenstein}. We record the following, which is a
consequence of \ref{EhregularMatlisCondition}.

\begin{prop}\label{SpecialGorensteinConditions}
  Suppose that the map $R\to k$ of \rings/ is \ehregular.
  Then $R\to k$ is \hGorenstein/ of shift $\SH$ if and only if
  $\HomR(k,R)$ is equivalent to $\Shift^\SH k$ as a left $k$-module.
\end{prop}

The rest of the section provides techniques for recognizing
\hGorenstein/ homomorphisms $R\to k$.

\begin{prop}\label{GorensteinMooreDuality}
  Suppose that $R$ is an \augmented/ $k$-algebra, and let
  $\EE=\EndR(k)$.  Assume that $\Pair Rk$ is \dccomplete, and that
  $R\to k$ is \ehregular. Then $R\to k$ is \hGorenstein/ (of shift
  $a$) if and only if
  $\EE\to k$ is \hGorenstein/ (of shift $a$). \qed 
\end{prop}

  See \cite[2.1]{halperin-thomas} for a differential graded version of this.

\begin{proof}
 Compute 
  \[
   \begin{aligned}
    \HomR(k,R)&\weq\HomR(k,\HomE(k,k))\weq
    \Hom_{R\Tensork\EE}(k\Tensork k,k)\\
    \HomE(k,\EE)&\weq\HomE(k,\HomR(k,k))\weq\Hom_{\EE\Tensork
      R}(k\Tensork k, k)\,.
    \end{aligned}
  \]
  What we are using is the
  fact that if $A$ and $B$ are left modules over the $k$-algebras $S$
  and $T$, respectively, and $C$ is a left module over $S\Tensor_kT$,
  then there is an equivalence
  \[ \Hom_S(A, \Hom_T(B, C)) \weq \Hom_{S\Tensork T}(A\Tensork B,C)\,.\]
  In our case $A=B=C=k$.
  This reveals a subtlety: $k\Tensork k$ is certainly equivalent to
  $k$, but not necessarily in a way which relates the tensor product
  action of $R\Tensork \EE$ on $k\Tensork k$ to the action of
  $R\Tensork \EE$ on $k$ given by $\EE=\EndR(k)$. Nevertheless, it is
  clear that $\HomR(k,R)$ is equivalent to a shift of $k$ if and only if
  $\HomE(k,\EE)$ is. If $\EE$ is \hGorenstein, then $R$ is
  \hGorenstein/ by \ref{SpecialGorensteinConditions}. If $R$ is
  \hGorenstein, $\EE$ is \hGorenstein/ by \ref{RegularCriterion} and
  \ref{SpecialGorensteinConditions}.
\end{proof}

\begin{prop}\label{HomInTheMiddle}
  Suppose that $S\to R$ is a map of \augmented/ commutative
  $k$-algebras such that $R$ is small as an $S$-module. Let
  $Q$ be the \augmented/ $k$-algebra $k\TensorOver S R$. Then
  there is an equivalence of $k$-modules
  \[
    \HomOver R(k,R)\weq\HomOver Q(k, \HomOver S(k,S)\Tensork Q),
  \]
  where the action of $Q$ on $\HomOver S(k,S)\Tensork Q$ is
  induced by the usual action of $Q$ on itself.
\end{prop}

There is a rational version in \cite[4.3]{halperin-thomas}. 
The argument below depends on the following general lemma, whose proof we leave
to the reader.

\begin{lem}\label{HochschildTensorHom}
  Suppose that $R$ is a $k$-algebra, that $A$ is a right $R$-module,
  and that $B$ and $C$ are left $R$-modules. Then there are natural
  equivalences
  \[
  \begin{aligned}
    \HomR(B,C)&\weq\HomOver{R\Tensork R\op}(R,\Homk(B,C))\\
    A\TensorR B &\weq R\TensorOver{R\Tensork R\op} (A\Tensork B)\,.
  \end{aligned}
  \]
\end{lem}

\begin{titled}{Proof of \ref{HomInTheMiddle}}
  Since $R$ is commutative, we do not distinguish in notation between
  $R$ and $R\op$.
  First note that
  \[
     \HomOver R(k,R)\weq \HomOver{R\TensorOver S R}(R, \HomOver S(k,R))
  \]
  as in \ref{HochschildTensorHom}. Now observe that $R$ is small over
  $S$, so that
  \begin{equation}\label{KeySmallnessAssumption}
     \HomOver S(k,R)\weq\HomOver S(k,S)\TensorOver S R\,.
  \end{equation}
  Under this equivalence, the left action of $R$ on
  $\HomOver S(k,S)\TensorOver S R$ is induced by the left action of $R$ on
  itself, and the right action of $R$ by the
  left action of $R$ on $k$.  Now since $S$ is commutative, the right and
  left actions of $S$ on $\HomOver S(k,S)$ are the same. In particular, the
  right action (which is used in forming $\HomOver S(k,S)\TensorOver S R$)
  factors through the homomorphism $S\to k$, and we obtain an
  equivalence
  \begin{equation}\label{KeySplittingFormula}
     \HomOver S(k,S)\TensorOver S R\weq \HomOver S(k,S) \Tensork
     (k\TensorOver S R)\weq \HomOver S(k,S)\Tensork Q\,.
  \end{equation}
  Let $M=\HomOver S(k,S)\Tensork Q$. Under \ref{KeySmallnessAssumption} and
  \ref{KeySplittingFormula} the left action of $R$ on $M$ is induced
  by the left action of $R$ on $Q$, while the right action of $R$ is
  induced by the left action of $R$ on $k$.  In particular, the
  action of $R\TensorOver S R$ on $M$ factors through an action of
  $k\TensorOver S R\weq Q$ on $M$, and so by adjointness we have
  \[
    \begin{aligned}
    \HomOver {R\TensorOver S R}(R,M) &\weq \HomOver{Q}
                (Q\TensorOver{R\TensorOver S R} R, M)\\
                    &\weq \HomOver Q (k, M) \,,
     \end{aligned}
   \]
   where the last equivalence depends on the calculation
   (\ref{HochschildTensorHom}) 
   \[
        (k\TensorOver S R)\TensorOver
        {R\TensorOver S R} R \weq k\TensorOver R R\weq k\,.
   \]
   The action of $Q$ on this object is the obvious one that factors
   through $Q\to k$. Combining the above gives the desired statement. \qed
\end{titled}

\begin{prop}\label{MixedFibration}
  Let $S\to R$ be a homomorphism of commutative \augmented/
  $k$-algebras, and set $Q=k\TensorOver S R$. Suppose that $R$ is
  small as an $S$-module, and that $R\to k$ is \ehregular. Then if the
  maps $S\to k$ and $Q\to k$ are \hGorenstein, so is $R\to k$.
\end{prop}

\begin{proof}
  By \ref{HomInTheMiddle}, $\HomOver R(k,R)\weq \Shift^\SH k$. It follows from
  \ref{SpecialGorensteinConditions} that $R\to k$ is \hGorenstein.
\end{proof}

\begin{NumberedSubSection}{Poincar\'e Duality}\label{DefPoincareDuality}
 
  A $k$-algebra $R$ is said to \emph{satisfy Poincar\'e duality} of
  dimension $\SH$ if there is an  equivalence
  $R\to\Shift^{\SH}\Homk(R,k)$ of $R$-modules; note that here we give $\Homk(R,k)$
  the left $R$-module structure induced by the right action of $R$ on
  itself. The algebra $R$ satisfies this condition if and only if
  there is an orientation class $\omega\in\Hof{-\SH}\Homk(R,k)$ with the
  property that $\Hof*\Homk(R,k)$ is a free module of rank one over
  $\Hof*R$ with generator $\omega$. If $k$ is a field, then
  $\Hof*\Homk(R,k)=\hom_k(\pi_*R,k)$, and $R$ satisfies Poincar\'e
  duality if and only if $\pi_*R$ satisfies Poincar\'e duality in the
  simplest algebraic sense.

  \begin{prop}\label{PoincareGivesGorenstein}
    Suppose that $R$ is \anaugmented/ $k$ algebra such that the map
    $R\to k$ is \ehregular. If $R$ satisfies Poincar\'e duality of
    dimension $\SH$, then $R\to k$ is \hGorenstein/ of shift~$\SH$. 
  \end{prop}

  \begin{proof}
    As in \ref{kAlgebraExample}, compute \[
          \HomR(k,R)\weq\HomR(k,\Shift^{\SH}\Homk(R,k))\weq
    \Shift^{\SH}\HomOver{k}(R\Tensor_Rk,k)\weq \Shift^{\SH}k\,.\] The fact
    that $R\to k$ is \hGorenstein/ follows from
    \ref{SpecialGorensteinConditions}.
  \end{proof}

We now give a version of the result from commutative ring theory that
``regular implies Gorenstein''.

  \begin{prop}\label{RegularGivesGorenstein}
    Suppose that $k$ is a field, $R$ is a connective commutative
    \ring, and $R\to k$ is \anhregular/ homomorphism which is surjective on
    $\Hof0$. Assume that the pair $\Pair Rk$ is \dccomplete. Then
    $R\to k$ is \hGorenstein/.
  \end{prop}

  \begin{rem}\label{SpecialRegularGivesGorenstein}
    It is possible to omit the \dccompleteness/ hypothesis from
    \ref{RegularGivesGorenstein} in the commutative ring case.
    Suppose that $R$ is a commutative Noetherian ring, $I\subset R$ is
    a maximal ideal, $k=R/I$ is the residue field, and
    $R\to k$ is \hregular. We show that $R\to k$ is also \hGorenstein. To see
    this, let $S=\lim_sR/I^s$ be the $I$-adic completion of $R$.
    As in the proof of \ref{ThreeCompletionTypes}, $S$ is flat
    over $R$ and $\Tor_0^R(S,k)\iso k$; in addition, the map
    $R\to S$ is a $k$-equivalence (of $R$-modules). This gives a
    chain of equivalences
    \[
    \HomR(k,R)\weq \HomR(k,S)\weq\HomOver{S}(S\TensorR
      k,S)\weq \HomOver{S}(k,S)\,.
    \]
    The flatness easily implies that $S\to k$ is also \hregular,
    and so $R\to k$ is \hGorenstein/ if and only if $S\to k$ is
    \hGorenstein. But it follows from \ref{ThreeCompletionTypes} that
    the pair 
    $\Pair Sk$ is \dccomplete, and so $S\to k$ is \hGorenstein/ by
    \ref{RegularGivesGorenstein}.
  \end{rem}

  \begin{rem}
    As the arguments below suggest, \ref{RegularGivesGorenstein} fails
    without the commutativity assumption. For instance, let $k=\Fp$, let $X$ be a
    simply-connected finite complex which does not satisfy Poincar\'e
    duality, and let $R\to k$ be the augmentation map $C_*(\Omega X;k)\to k$.
    Then $R\to k$ is \hregular/ (\ref{FiniteSkeletaFiniteStableType}),
    $(R,k)$ is \dccomplete/ (\ref{ChainCochainStuff}), but it is easy
    to see that $R\to k$ is not \hGorenstein.
  \end{rem}

  \begin{lem}\label{AdamsLemma}
    Suppose that $k$ is a field, $R$ is a connective commutative
    \ring, and $R\to k$ is a homomorphism which is
    surjective on $\Hof0$. Assume that $k$ is of upward finite type
    over $R$. Then $\Hof*\EndR(k)$ is in a natural way a cocommutative
    Hopf algebra over $k$.
  \end{lem}

  \begin{proof}
    The diagram chasing necessary to prove this is described in detail
    in \cite[pp.~56--76]{adams-lectures}, with a focus at the end on
    the case in which $R=\Sphere$, $k=\Fp$, and $\Hof *\EndR(k)$ is the
    mod~$p$ Steenrod algebra. Let $\EE=\EndR(k)$. The key idea is that
    $\Hof*\EE$ is the $k$-dual of the commutative $k$-algebra
    $\Hof*(k\TensorR k)$: as in \ref{kAlgebraExample} there are
    equivalences
    \[
      \Homk(k\TensorR k,k) \weq \HomR(k, \Homk(k,k))\weq\EndR(k)\,\,.
    \]
    The $k$-dual of the multiplication on $\Hof*(k\TensorR
    k)$ then provides the comultiplication on $\Hof*\EndR(k)$. The
    fact that $k$ is of upward finite type over $R$ guarantees that
    the groups $\Hof i(k\TensorR k)$ are finite-dimensional over $k$.

    There is a technicality: $k\TensorR k$ is a bimodule over $k$, not an
    algebra over $k$. However $k\TensorR k$ \emph{is} an algebra over
    $R$, so that the surjection $\Hof0R\to k$ guarantees that the left
    and right action of $k$ on $\Hof*(k\TensorR k)$ agree. For the
    same reason, the left and right actions of $k$ on $\Hof*\EndR(k)$
    agree, and this graded ring becomes a Hopf algebra over $k$.
  \end{proof}

  \begin{titled}{Proof of \ref{RegularGivesGorenstein}}
    Let $\EE=\EndR(k)$. The connectivity assumptions on $R$ imply that
    $\Hof0\EE\iso k$ and that $\EE$ is coconnective; by
    \ref{AdamsLemma}, $\EE$ is a Hopf algebra over $k$.  In fact, $\EE$
    is finitely built from $k$ (\ref{RegularCriterion}), and so
    $\Hof*\EE$ is a finite dimensional Hopf algebra over $k$.
    Sweedler has remarked that a connected finite-dimensional Hopf
    algebra over $k$ with commutative comultiplication and
    involution satisfies algebraic Poincar\'e duality
    \cite{MoorePeterson}; see also
    \cite[5.1.6]{sweedler}.
    The map $\EE\to k$ is thus \hGorenstein/
    by \ref{PoincareGivesGorenstein}, and $R\to k$ by
    \ref{GorensteinMooreDuality}.  \qed
  \end{titled}

  \begin{rem}
    The above arguments are related to those of Avramov and Golod
    \cite{avramov-golod}, who show that a commutative Noetherian local ring $R$ is
    Gorenstein if and only if the homology of the associated Koszul
    complex is a Poincar\'e duality algebra.
  \end{rem}

\end{NumberedSubSection}

\section{A local cohomology theorem}\label{LocalCohomologyTheorem}

One of the attractions of the \hGorenstein/ condition on \aring/ $R$
is that it has structural implications for $\Hof *R$, which can
sometimes be thought of as a duality property. To illustrate this, we
look at the special case in which $R\to k$ is \anhGorenstein/ map of
\augmented/ $k$-algebras, where $k$ is a field. Let $\EE=\EndR(k)$. By
\ref{GorensteinMeansCells}, the \hGorenstein/ condition gives
\[
   \Shift^\SH k\TensorE k\weq \CellOf kR\,.
\]
We next assume that the
right $\EE$-structure on $\Shift^\SH k$ given by $\Shift^\SH k\weq\HomR(k,R)$ is
equivalent to the right $\EE$-structure given by \[\Shift^\SH
k\weq\HomR(k,\Shift^\SH \Homk(R,k))\,.\]
By \ref{MatlisTypeCellularApprox} this gives an equivalence
\[
   \Shift^\SH k\TensorE k \weq \Shift^\SH\CellOf k\Homk(R,k)\,.
\]
Assume in addition that $\Homk(R,k)$ is itself $k$-cellular as an
$R$-module. Combining the above then gives
\begin{equation}\label{CellsVsDuals}
   \Shift^\SH\Homk(R,k)\weq \CellOf kR\,.
\end{equation}
Now in some reasonable circumstances we might expect a spectral
sequence
\begin{equation}\label{MythicalSpectralSequence}
    E_{i,j}^2=\Hof i \CellOverOf {\Hof *R}k(\Hof*R)_j\Rightarrow\Hof{i+j}\CellOverOf Rk(R)
\end{equation}
which in the special situation we are considering would give
\[
E^2_{i,j}=\Hof i\CellOverOf {\Hof*R}k(\Hof*R)_j\Rightarrow
\Hof{i+j-\SH}\Homk(R,k)\,.
\]
(The subscript $j$ refers to the $j$'th homogeneous component of an
appropriate grading on $\Hof i\CellOverOf{\Hof*R}k(\Hof*R)$.) This is
what we mean by a duality property for $\Hof*R$: a spectral sequence starting from
some 
covariant algebraic data associated  to $\Hof*R$ and abutting to the
dual object 
$\Hof*\Homk(R,k)\iso\HomOver k(\Hof*R,k)$. If $R$ is $k$-cellular as a
module over itself, then \ref{CellsVsDuals} gives
$\Shift^\SH\Homk(R,k)\weq R$, and we obtain ordinary Poincar\'e
duality.

The problematic point here is the existence of the spectral sequence
\ref{MythicalSpectralSequence}. Rather than trying to construct this
spectral sequence in general and study its convergence properties, we
concentrate on a special case in which it is possible to identify
$\CellOverOf Rk(R)$ explicitly. To connect the following statement
with \ref{MythicalSpectralSequence}, recall
\cite[\S6]{dwyer-greenlees} that if $S$ is a commutative ring and
$I\subset S$ a finitely generated ideal with quotient ring $k=S/I$,
then for any discrete $S$-module $M$ the local cohomology group
$H_I^i(M)$ can be identified with $\Hof{-i}\CellOverOf Sk(M)$.

\begin{prop}\label{LocalCohomologySS}
  Suppose that $k$ is a field, and that $R$ is a coconnective
  commutative augmented $k$-algebra. Assume that $\Hof*R$ is
  Noetherian, and that the augmentation map induces an isomorphism
  $\Hof0R\iso k$. Then for any $R$-module $M$ there is a spectral
  sequence 
  \[
     E^2_{i,j}=H_I^{-i}(\pi_*M)_j\Rightarrow\Hof {i+j}\CellOverOf Rk(M)\,.
  \]
\end{prop}

Under the assumptions of \ref{LocalCohomologySS}, $\HomR(R,k)$ is
$k$-cellular as an $R$-module (\ref{BuiltFromk}).
Given the above discussion, this leads to the following
local cohomology theorem.

\begin{prop}\label{GoodLCSS}
  In the situation of \ref{LocalCohomologySS}, assume in addition that
  $R\to k$ is \hGorenstein/ of shift $\SH$, and that $k$ has a unique
  \Elift/ (where $\EE=\EndR(k)$).
  Then there is a spectral sequence
  \[
    E^2_{i,j}=H_I^{-i}(\Hof*R)_j\Rightarrow\Hof{i+j-\SH}\Homk(R,k)\,.
  \]  
\end{prop}

\begin{rem}
  The structural implications of this spectral sequence for the
  geometry of the ring $\Hof*R$ are investigated in
  \cite{greenlees-lyubeznik}. 
\end{rem}

\begin{titled}{Proof of \ref{LocalCohomologySS}}
We first copy some constructions from \cite[\S6]{dwyer-greenlees}.
For any $x\in\Hof{|x|}R$ we can form an $R$-module $R[1/x]$ by taking the
homotopy colimit of the sequence
\[
   R\RightArrow x \Shift^{-|x|}R\RightArrow x \Shift^{-2|x|}R\RightArrow x\cdots\,.
\]
(Actually, $R[1/x]$ can also be given the structure of a commutative
\ring, in such a way that $R\to R[1/x]$ is a homomorphism.) Write
$K_m(x)$ for the fibre of $x^m:R\to \Shift^{-m|x|}R$, and $K_\infty(x)$ for the
fibre of the map $R\to R[1/x]$. Now choose a finite sequence $
x_1,\ldots,x_n$ of generators for $I\subset\Hof*R$, and let
\[
\begin{aligned}
K_m&=K_m(x_1)\TensorR\cdots\TensorR K_m(x_n)\\
K_\infty&=K_\infty(x_1)\TensorR\cdots\TensorR K_\infty(x_n)\,.
\end{aligned}
\]
Recall that $R$ is commutative, so that right and left $R$-module
structures are interchangeable, and tensoring two $R$-modules over $R$
produces a third $R$-module. 
Write $K=K_1$. It is easy to see that $\Hof*K$ is finitely built from $k$ as a
module over $\Hof *R$, and hence (\ref{BestOfAllCoregularity}) that
$K$ is finitely built from $k$ as a module over $R$. An inductive
argument (using cofibration sequences $K_m(x_i)\to K_{m+1}(x_i)\to
\Shift^{m|x_i|}K_1(x_i)$) shows that $K$ builds $K_m$ and hence also builds
$K_\infty\weq\hocolim K_m$ (cf. \cite[6.6]{dwyer-greenlees}). It is
easy to see that the evident map $K_\infty\to R$ gives equivalences
\begin{equation}\label{InfiniteKoszulProperty}
k\TensorR K_\infty\weq k\quad\quad K\TensorR K_\infty\weq K\,.
\end{equation}
See \cite[proof of 6.9]{dwyer-greenlees}; the second equivalence
follows from the first because $K$ is built from $k$.  The first
equivalence implies that $K_\infty$ builds $k$ and this in turn shows that
$K$ builds $k$. Since $K$ is small over $R$, we see that $R\to k$ is
\ehregular/ with Koszul complex $K$. In particular, a map $A\to B$
of $R$-modules is a $k$-equivalence if and only if it is a
$K$-equivalence, or (since $\HomR(K(x_i), R)\weq\Shift K(x_i)$ and
hence $\HomR(K,R)\weq\Shift^nK)$) if and only if it induces an
equivalence $K\TensorR A\to K\TensorR B$. Since $K_\infty$ is built
{}from $k$ as an $R$-module, so is $K_\infty\TensorR M$.  The right hand
equivalence in \ref{InfiniteKoszulProperty} implies that the map
$K_\infty\TensorR M\to M$ induces an equivalence
\[
    K\TensorR K_\infty\TensorR M \to K\TensorR M\,\,,
\]
and it follows that $K_\infty\TensorR M$ is $\CellOverOf Rk(M)$.
Each module $K_\infty(x_i)$ lies in a cofibration sequence
\[
   \Shift^{-1}R[1/x_i]\to K_\infty(x_i)\to R
\]
which can be interpreted as a one-step increasing filtration of
$K_\infty(x_i)$. Tensoring these together gives an $n$-step filtration
of $K_\infty$,
\[
    0=F_{n+1}\to F_n\to F_{n-1}\to\cdots\to F_0=K_\infty,
\]
with the property that there are equivalences
\[
    F_s/F_{s+1} \weq  \bigoplus_{\{i_1,\ldots,i_s\}}
    R[1/x_{i_1}]\TensorR\cdots\TensorR R[1/x_{i_s}] \,.
\]
The sums here are indexed over subsets of cardinality $s$ from
$\{1,\ldots,n\}$. Tensoring this filtration with $M$ gives a
finite filtration of $\CellOverOf Rk(M)$, and the spectral sequence of
the proposition is the homotopy spectral sequence associated to the
filtration. The identification of the $E^2$-page as local cohomology
is standard \cite[\S6]{dwyer-greenlees} \cite{GrothendieckLocal}; the
main point here is to notice that since $\Hof *R[1/x_i]$ is flat over
$\Hof*R$, there are isomorphisms $\Hof*(R[1/x_i]\TensorR
M)\iso\Hof*(R[1/x_i])\TensorOver{\Hof*R}\Hof*M\iso
(\Hof*R)[1/x_i]\TensorOver{\Hof*R}\Hof*M$. \qed
\end{titled}

\section{Gorenstein examples}\label{GorensteinExamples}

We give several examples of \rings/ which are \hGorenstein, and at
least one example of \aring/ which is not. Of course,
commutative Noetherian local rings which are not Gorenstein are easy
to come by; these also provide examples of non--Gorenstein \rings.

\begin{NumberedSubSection}{\Hregular/ chains}\label{GoodFiniteComplex}
  Suppose that $X$ is a pointed connected topological space and that
  $k$ is a field such that the pair $(X,k)$ is of Eilenberg-Moore type
  (\ref{ChainCochainStuff}), and such that  $H^*(X,k)$ is
  finite--dimensional and satisfies Poincar\'e duality of formal
  dimension~$b$. Let $a=-b$. We point out that
  the augmentation map $C_*(\Omega X;k)\to k$ is \hregular/ and
  \hGorenstein/ of shift~$a$; note that the shift~$a$ is negative in
  this case.

  Let $R=C^*(X;k)$, $\EE=C_*(\Omega
  X;k)\weq\EndR(k)$.
  By \ref{RegularCochainAlgebras}, $R\to k$ is \cohregular/ and $\EE\to k$ is
  \hregular/. The map  $R\to k$ is
  \hGorenstein/ of shift $\SH$ (\ref{PoincareGivesGorenstein}) and
  $\EE\to k$ is also \hGorenstein/ with the same shift
  (\ref{GorensteinMooreDuality}).  The ring $R$ has a local cohomology
  spectral sequence (\ref{GoodLCSS}), but this collapses to a
  restatement of Poincar\'e duality:
  \[
     E^2=\Hof*R\iso\CellOverOf {\Hof*R} k(\Hof*R)\iso
     \Hof*\Sigma^{a}\Homk(R,k)\,.
   \]
   In the absence of the hypothetical spectral sequence
   \ref{MythicalSpectralSequence}, there is nothing like a local
   cohomology theorem for the noncommutative \ring/ $\EE$.
\end{NumberedSubSection}

\begin{NumberedSubSection}{\Hregular/ cochains}\label{GoodLoopSpace} 
  Suppose that $k$ is a field and $G$ is a topological group such that
  $H_*(G;k)$ is finite dimensional. Assume in addition that $(BG,k)$ is
  of Eilenberg-Moore type; this covers the cases in which $k=\Fp$, and
  $G$ is a finite $p$-group, a compact Lie group with $\pi_0G$ a
  finite $p$-group, or a $p$-compact group
  \cite{dwyer-wilkerson;methods}. We point out that the augmentation
  map $C^*(BG;k)\to k$ is \hregular/ and \hGorenstein, and satisfies a
  local cohomology theorem.

Let $R=C^*(BG;k)$ and
  $\EE=C_*(G;k)$. The map $\EE\to k$ is
  \cohregular/ (\ref{BestOfAllCoregularity}), and hence $R\to k$ is \hregular/
  (\ref{RegularChainAlgebras}).  The graded ring $H_*(G;k)$ is a
  finite dimensional group-like Hopf algebra over $k$, and so by
  Sweedler (cf. \cite[5.1.6]{sweedler}) satisfies algebraic Poincar\'e
  duality of some dimension, say $\SH$. If $G$ is a connected compact
  Lie group, then $\SH=\dim G$; the ``fundamental class'' $\omega$
  (see \ref{DefPoincareDuality}) lies in \[H^{\SH}(G;k)=\pi_{-\SH}C^*(G;k)=\Hof{-\SH}(k\Tensor_Rk)\,.\] By
  \ref{PoincareGivesGorenstein}, $\EE\to k$ is \hGorenstein/ of shift
  $\SH$, and so $R\to k$ is also \hGorenstein/ with the same shift
  (\ref{GorensteinMooreDuality}).  The graded ring
  $H^*(BG;k)=\Hof*R$ is Noetherian. If $k$ is of characteristic
  zero, this follows from the fact that the ring is a finitely
  generated polynomial algebra over $k$; see
  \cite[7.20]{milnor-moore}. If $k=\Fp$ and $G$ is a compact Lie
  group, the finite generation statement is a classical theorem of
  Golod \cite{golod} and Venkov \cite{venkov}; in
  the $p$-compact group case it amounts to the main result of
  \cite{dwyer-wilkerson;methods}. By \ref{BuiltFromk} and \ref{GoodLCSS} there is a local
  cohomology theorem for $R$.
\end{NumberedSubSection}

\begin{NumberedSubSection}{Completed classifying spaces}\label{GorensteinCompactLie}
  Suppose that $G$ is a  compact Lie group such that the adjoint
  action of $G$ on its Lie algebra is orientable (e.g., 
  $G$ might be a finite group). Let $k=\Fp$.  
   We point out that the augmentation map $C^*(BG;k)\to k$ is
  \ehregular/ and \hGorenstein, and that $C^*(BG;k)$ has a local
  cohomology theorem.   Note that we are not assuming 
  that $(BG,k)$ is of Eilenberg-Moore
  type, and in particular we are not assuming that $\pi_0G$ is a
  $p$-group; the present case contrasts with \ref{GoodLoopSpace} in that
  $C^*(BG;k)\to k$ need not be \hregular.

   We continue the discussion in
  \ref{EHRegularChains}, with the same notation.  Recall that $X$ is
  the $p$-completion of $BG$, $R=C^*(X;k)\weq C^*(BG;k)$, and
  $\EE=C_*(\Omega X; k)$; the space $\Omega X$ plays the role of $G$
  above in \ref{GoodLoopSpace}, but we do not have that $H_*(\Omega
  X;k)$ is finite dimensional. The fibre
  $M$ in \ref{EmbeddingFibration} is a compact manifold;
  it is orientable because its tangent bundle is the bundle associated
  to the conjugation action of $G$ on the quotient of the  Lie algebra of
  $SU(n)$ by the Lie algebra of $G$. (Note that 
  since $SU(n)$ is connected, the conjugation action of $G$ on the Lie
  algebra of $SU(n)$ preserves
  orientation.) As in \ref{GoodFiniteComplex}, $Q=C^*(M;k)$ is
  \cohregular/ and \hGorenstein. Similarly, $S=C^*(BSU(n);k)$ is
  \hregular/ and \hGorenstein/ by \ref{GoodLoopSpace}. Let
  $S'=C_*(SU(n);k)$, so that $S\weq\End_{S'}(k)$. The group $SU(n)$
acts on $M=SU(n)/G$, so that $S'$ acts on $C_*(M;k)$; by an
Eilenberg-Moore spectral sequence argument, there is an equivalence
$R\weq\Hom_{S'}(C_*(M;k),k)$. It follows from
  \ref{DownwardSmall} that $R$ is small as a module over $S$. By
  \ref{EhregularInTheMiddle} and \ref{MixedFibration}, $R\to k$ is \ehregular/ and \hGorenstein/, as
  is $\EE\to k$ (\ref{GorensteinMooreDuality}). Since $\Homk(R,k)$ is
  $k$-cellular over $R$ (\ref{BuiltFromk}), there is a local
  cohomology spectral sequence for $R$ (\ref{GoodLCSS}).
\end{NumberedSubSection}

\begin{NumberedSubSection}{Finite complexes}\label{FiniteComplexes}
  Suppose that $X$ is a pointed connected finite complex which is a
  Poincar\'e duality complex over $k$ of formal dimension $b$; in
  other words, assume that $X$ satisfies possibly unoriented Poincar\'e
  duality with arbitrary (twisted) $k$-module coefficients. To be
  specific, assume that $k$ is a finite field, the field $\Q$, or the
  ring $\Z$ of integers. We point out that the augmentation map
  $C_*(\Omega X;k)\to k$ is \hregular/ and \hGorenstein. This is
  related to \ref{GoodFiniteComplex} but slightly different: here we
  assume that $X$ is finite, we do not insist that $(X,k)$ be of
  Eilenberg-Moore type, but we require (possibly twisted) Poincar\'e
  duality with arbitrary $k$-module coefficients.

  Let $R$ denote the augmented $k$-algebra
  $C_*(\Omega X;k)$, so that $\Hof0R\iso k[\pi_1X]$. Note that $R\to
  k$ is \hregular/ (\ref{FiniteSkeletaFiniteType}). Any module $M$
  over $k[\pi_1X]$ gives a module over $R$, and (by a version of the
  Rothenberg-Steenrod construction) there are isomorphisms
  \[
     H_i(X;M)\iso \Hof i (k\TensorR M) \quad\quad
     H^i(X;M)\iso\Hof{-i}\HomR(k,M)\,\,.
  \]
  Let $a=-b$. 
  The duality condition on $X$ can be expressed by saying that there
  is a module $\lambda$ over $k[\pi_1X]$ whose underlying $k$-module
  is isomorphic to $k$ itself, and an orientation class $\omega\in
  \Hof{-\SH} (\lambda\TensorR k)$, such that, for any $k[\pi_1X]$-module $M$,
  evaluation on $\omega$ gives an equivalence
  \begin{equation}\label{TwistedDualityEquiv}
        \HomR(k, M) \to \Shift^{\SH} \lambda\TensorR M\,\,.
  \end{equation}
  By \ref{UniqueModuleStructures}, it follows 
  that \ref{TwistedDualityEquiv} is an equivalence for any $R$-module
  $M$ which has only one nonvanishing homotopy group. By triangle
  arguments (cf. \ref{ConnectivePostnikov}) it is easy to conclude
  that \ref{TwistedDualityEquiv} is an equivalence for all $M$ which
  have only a finite number of nonvanishing homotopy groups, and by
  passing to a limit (cf. proof of \ref{PChainCanonical}) that
  \ref{TwistedDualityEquiv} is actually an equivalence for all
  $R$-modules $M$. Note that this passage to the limit depends on the
  fact that $k$ is small over $R$. The case $M=R$ of
  \ref{TwistedDualityEquiv} gives
  \[\HomR(k,R)\weq \Shift^{\SH}\lambda\TensorR R\weq\Shift^{\SH}\lambda\weq
  \Shift^{\SH}k\,\,,\] and so by \ref{SpecialGorensteinConditions}, $R\to k$ is
  \hGorenstein/ of shift $\SH$. Let $\EE=\EndR(k)$. The pair $(R,k)$
  is not necessarily \dccomplete, and so $\EE\to k$ is not necessarily
  \hGorenstein; for example, it is clear that $\Hof*\EE\iso H^*(X;k)$
  need not satisfy algebraic Poincar\'e duality in the nonorientable case.

  The equivalence $\HomR(k,R)\weq\lambda$ is an $R$-module
  equivalence as long as $\HomR(k,R)$ is given the right $R$-module structure
  obtained from the right action of $R$ on itself. In this way the
  orientation character of the Poincar\'e complex $X$ is derived from
  the one bit of structure on $\HomR(k,R)$ that does not play a role
  in the definition of what it means for $R\to k$ to be \hGorenstein/
  (\ref{GorensteinMeansNiceDuals}).
\end{NumberedSubSection}

\begin{NumberedSubSection}{Suspension spectra of loop
    spaces}\label{SuspensionSpectra}
  This is a continuation of \ref{FiniteComplexes} and \ref{MatlisSuspensionLoops}:
 $X$ is a pointed connected finite complex which is a
  Poincar\'e duality complex over $k$ of formal dimension $\SH$. Let $k=\Sphere$.
  We observe that the augmetation map $C_*(\Omega X;k)\to k$ is
  \hregular/ and \hGorenstein, and point out how this Gorenstein
  condition leads to a homotopical construction of the Spivak normal
  bundle of~$X$.

   Let $R$ denote the
  augmented $k$-algebra $C_*(\Omega X;k)$, and let $\EE$ denote
  $C^*(X;k)\weq\EndR(k)$.  The map $R\to k$ is \hregular/ by \ref{FiniteSkeletaFiniteType}.
    Note that
  $S=\Z\TensorS R\weq C_*(\Omega X;\Z)$ (\ref{ChainCochainStuff}).
  We wish to show that $R\to k$ is \hGorenstein/ of shift
  $-\SH$, or equivalently, that $\HomR(k,R)\weq \Shift^{-\SH} k$. The
  spectrum $Y=\Shift^{-\SH}k$ is characterized by a combination of the
  homotopical property that $Y$ is bounded below, and the homological
  property that $\Z\Tensork Y\weq\Shift^{-\SH}\Z$. The spectrum
  $\HomR(k,R)$ is bounded below because $R$ is bounded below and $k$
  is small over $R$. Similarly, the fact that $k$ is small over $R$
  implies that $\Z\Tensork\HomR(k,R)\weq\HomR(k, \Z\Tensork R)$. Now
  compute
  \[
    \HomR(k,\Z\Tensork R) \weq \HomOver{\Z\Tensork R}(\Z, \Z\Tensork
    R)\weq \Shift^{-\SH}\Z\,\,,
  \]
  where the first equivalence comes from adjointness, and the second
  from \ref{FiniteComplexes}. It follows that $R\to k$ is \hGorenstein. (If
  $X$ is simply connected, then $\EE\to k$ is \cohregular/ and
  \hGorenstein, since in this case $(R,k)$ is \dccomplete.)

  The stable homotopy orientation character of $X$ is given by the
  action of $R$ on $k\weq\Sphere$ obtained via
  $\Shift^{-\SH}k\weq\HomR(k,R)$ from the right action of $k$ on
  itself; see \ref{TwistedDualityEquiv} for the homological version of
  this. It is not too far off to interpret this character as a
  homomorphism $\Omega X\to \Sphere^\times$; in any case it determines
  a stable spherical fibration over $X$ which can be identified with
  the Spivak normal bundle. (To see this, note that the Thom complex of this spherical
  fibration is $\HomR(k,R)\TensorR k\weq\HomR(k,k)=\EE$, and the top
  cell has a spherical reduction given by the unit homomorphism
  $\Sphere\to\EE$.) For some more details see \cite{klein}.
\end{NumberedSubSection}

\begin{NumberedSubSection}{The sphere spectrum}
  Let $R=\Sphere$ and $k=\Fp$. The map $R\to k$ is not \hGorenstein; in
  fact, by Lin's theorem \cite{lin}, $\HomR(k,R)$ is trivial. Neeman
  has suggested that \Sphere/ should be considered to be a ``fairly
  ordinary non-Noetherian ring'' \cite[\S0]{rNeeman}. If so, then Lin's theorem
  is perhaps analogous to the classical calculation that if
  $R=\Fp[x_1,x_2,\ldots]$ and $k=R/I$ for 
  $I=\langle x_1,x_2,\ldots\rangle$, then $\HomR(k,R)=0$.
\end{NumberedSubSection}

\providecommand{\bysame}{\leavevmode\hbox to3em{\hrulefill}\thinspace}

\end{document}